\documentclass[aop,dvips]{arximspdf}
\usepackage{mathbh}
\usepackage{graphicx}

\doi{10.1214/009117906000000881}
\volume{35}
\issue{3}
\pubyear{2007}
\firstpage{833}
\lastpage{866}

\makeatletter
\newtheorem{prop}{Proposition}[section]
\newtheorem{theorem}[prop]{Theorem}
\newtheorem{coro}[prop]{Corollary}
\newtheorem{lemme}[prop]{Lemma}

\newcommand{\boite}{R}
\newcommand{\voisinage}{\mathcal{R}}

\newcommand{\communique}{\leftrightarrow}
\newcommand{\communiqueM}{\leftrightarrow}
\newcommand{\N}{\mathbb{Z}_{+}}

\newcommand{\Z}{\mathbb{Z}}

\newcommand{\Zd}{\mathbb{Z}^d}

\newcommand{\1}{\mathbh{1}}
\newcommand{\Q}{\mathbb{Q}}

\newcommand{\Pcond}{\overline{\mathbb{P}}}
\newcommand{\R}{\mathbb{R}}
\newcommand{\Rd}{\mathbb{R}^d}

\renewcommand{\P}{\mathbb{P}}

\newcommand{\Ed}{\mathbb{E}^d}

\newcommand{\Ber}{\mathrm{Ber}}

\newcommand{\Sd}{\ensuremath{\mathfrak{S}_{d}}}
\newcommand{\frinstar}{\partial_*^\mathrm{in}}
\newcommand{\as}{\mbox{ a.s.}}

\renewcommand{\epsilon}{\varepsilon}

\makeatother

\begin{document}
\begin{frontmatter}

\title{Large deviations for the chemical distance in supercritical Bernoulli percolation}
\runtitle{Large deviations for the chemical distance}

\begin{aug}
\author[A]{\fnms{Olivier}
\snm{Garet}\corref{}\ead[label=e1]{Olivier.Garet@univ-orleans.fr}\ead[url,label=u1]{http://www.univ-orleans.fr/mapmo/membres/garet/}} and
\author[B]{\fnms{R\'{e}gine} \snm{Marchand}\ead[label=e2]{Regine.Marchand@iecn.u-nancy.fr}
\ead[url,label=u2]{http://www.iecn.u-nancy.fr/\texttildelow marchand/}}
\runauthor{O. Garet and R. Marchand}
\affiliation{University of Orl\'eans and University of Nancy}
\address[A]{Laboratoire de Math{\'e}matiques, Applications \\
\quad et Physique
Math{\'e}matique d'Orl{\'e}ans \\
UMR 6628 Universit{\'e} d'Orl{\'e}ans\\
BP 6759\\
 45067 Orl{\'e}ans Cedex 2\\ France\\
\printead{e1}\\
\printead{u1}} 
\address[B]{Institut Elie Cartan Nancy\\\quad (math{\'e}matiques)\\
Universit{\'e} Henri Poincar{\'e} Nancy 1\\
Campus Scientifique\\ BP 239 \\
54506 Vandoeuvre-l{\`e}s-Nancy  Cedex\\ France\\
\printead{e2}\\
\printead{u2}}
\end{aug}

\received{\smonth{10} \syear{2005}}
\revised{\smonth{9} \syear{2006}}

\begin{abstract}
The chemical distance $D(x,y)$ is the length of the shortest open path
between two points $x$ and $y$ in an infinite Bernoulli percolation cluster.
In this work, we study the asymptotic behavior of this random metric,
and we prove that, for an appropriate norm $\mu$ depending on the dimension
and the percolation parameter,
the probability of the event
\[
\biggl\{\   0 \communique x, \frac{D(0,x)}{\mu(x)}\notin (1-\epsilon,
  1+\epsilon) \
\biggr\}
\]
exponentially decreases when $\|x\|_1$ tends to infinity. From this bound
we also derive a large deviation inequality for the corresponding asymptotic
shape result.
\end{abstract}


\end{frontmatter}

\section{Introduction and statement of main results}
The  matter of this article is the study of the asymptotic length of the
shortest open path between two points in an infinite Bernoulli percolation
cluster.

Let us first recall the Bernoulli percolation model and its usual notation.
Consider the
graph whose vertices are the points of
$\Z^d$, and put a nonoriented edge between each pair $\{x,y\}$ of points in $\Z^d$ such
that the Euclidean distance between~$x$ and $y$ is equal to $1$: two such
points are called \textit{neighbors}, and this set
of edges is denoted by $\Ed$.

Set $\Omega=\{0,1\}^{\Ed}$. We denote by $\P_p$ the product probability
$(p\delta_1+(1-p)\delta_0)^{\otimes \Ed}$ on the set $\Omega$.
For a  point $\omega$ in $\Omega$, we say that the edge $e \in \Ed$ is
\textit{open} in the configuration $\omega$ if $\omega(e)=1$, and \textit{closed}
otherwise. The states of the different edges are thus
independent under $\P_p$.
In the whole paper, the parameter $p$ is supposed to satisfy $p\in (p_c,1],$
where $p_c=p_c(d)$ is the critical probability for Bernoulli bond
percolation on $\Z^d$.

A \textit{path} is a sequence $\gamma=(x_1,
e_1,x_2,e_2,\ldots,x_n,e_n,x_{n+1})$ such that $x_i$ and $x_{i+1}$ are
neighbors and $e_i$ is the edge between $x_i$ and $x_{i+1}$.
We will also sometimes describe $\gamma$ only by the vertices it
visits $\gamma=(x_1,x_2,\ldots,x_n,x_{n+1})$ or by its edges
$\gamma=(e_1,e_2,\ldots,e_n)$. The number $n$ of edges in $\gamma$ is called
the \textit{length} of $\gamma$ and is denoted by $|\gamma|$. Moreover, we
will only consider \textit{simple paths} for which the visited vertices are all
distinct. A path is said to be \textit{open} in the configuration $\omega$
if all its edges are open in $\omega$.

The \textit{clusters} of a configuration $\omega$ are the connected
components of the graph induced on $\Z^d$ by the open edges in
$\omega$. For $x$ in $\Z^d$, we denote by $C(x)$ the cluster containing
$x$. In other words, $C(x)$ is the set of points in $\Z^d$ that are linked
to $x$ by an open path.
We write $x\communique y$ to signify that $x$ and $y$ belong to the same cluster.
For $p> p_c$, there exists almost surely one
and only one infinite cluster. We denote by $C_{\infty}$ the random set:
$C_{\infty}=\{k\in \Zd\dvtx  |{C(k)}|=+\infty\}$, which is
almost surely connected.

We introduce the \textit{chemical distance} $D(x,y)(\omega)$ between $x$ and
$y$ in $\Zd$,
depending on the Bernoulli percolation configuration $\omega$:
\[
D(x,y)(\omega)= \inf_{\gamma}|\gamma|,
\]
where the infimum is taken on the set of paths whose ends are $x$
and $y$ and that are open in the configuration $\omega$.
It is of course only defined when $x$ and $y$ are in the same percolation
cluster. Otherwise, we set by convention $D(x,y)=+\infty$.
The random distance $D(x,y)$ is thus, when it is finite, the minimal number of
open edges needed to link $x$ and $y$ in the configuration $\omega$, and is
thus larger than $\|x-y\|_1$, where $\|\cdot\|_1$ is the usual $\ell^1$ norm:
$\|x\|_1=\sum_{i=1}^d|x_i|$.

Note that the chemical distance $D(0,x)$ on the infinite Bernoulli cluster with
parameter $p>p_c$ can be seen as the travel time
between $0$ and $x$ in a first-passage percolation model where the passage times
of the edges  are independent identically distributed random variables with common distribution $p\delta_1+(1-p)\delta_{+\infty}$.

Antal and Pisztora~\cite{ap96} have proved that the chemical distance cannot
asymptotically be
too large when compared with the usual distance  $\|\cdot\|_1$: for each $p>p_c$,
there exists a positive constant $\rho$ such that
\begin{equation}
\label{AP}
\mathop{\lim\sup}_{\Vert x\Vert_1\to +\infty} \frac{\ln \P_p(0\communique
  x,   D(0,x)\ge \rho\Vert x\Vert_1)}{\Vert x\Vert_1}<0.
\end{equation}


If we think of the chemical distance as a special travel time in a
first-passage percolation model,
it is natural to expect that the term $\rho\|x\|_1$ in~(\ref{AP})
could be replaced by a smaller term, depending on a directional functional.
Indeed, a good candidate exists, which has been defined in a previous paper
of the authors:
\begin{prop}[(\cite{gm04})]
\label{puiskilfo}
Let $p>p_c$ and consider the chemical distan\-ce~$D(\cdot,\cdot)$ for Bernoulli percolation with parameter $p$.
There exists a norm $\mu$ on~$\Rd$ such that, almost surely on the event
$\{0\communique\infty\}$,
\[
\forall u\in\Zd\qquad\lim_{n \rightarrow \infty} \frac{D(0,T_{n,u}u)}{T_{n,u}}=\mu(u),
\]
where $(T_{n,u})_{n\ge 1}$ is the
increasing sequence of positive integers $k$
such that \mbox{$ku\communique\infty$}.
\end{prop}

It is then natural to study the fluctuations around this limit, and to look for exponential decay results analogous to the ones
obtained, for instance by Grimmett and Kesten~\cite{gk84}, in
first-passage percolation.
Therefore, the main objective of this work is to prove the following large
deviation bound:
\begin{theorem}
\label{ap-enrhume}
Let $p$ in the interval $(p_c,1]$ and denote by $\mu$
the norm on $\Rd$ given by Proposition~\textup{\ref{puiskilfo}}. Then,
\[
\forall \varepsilon>0 \qquad \mathop{\lim\sup}_{\Vert x\Vert_1\to +\infty} \frac{\ln \P_p  (0\communique
  x,   {D(0,x)}/{\mu(x)}\notin (1-\epsilon, 1+\epsilon) )}{\Vert x\Vert_1}<0.
  \]
\end{theorem}

The proof of Theorem~\ref{ap-enrhume} is divided into
two parts: the upper large deviations and the lower large deviations, which
are, respectively, dealt with in Sections \ref{upper}~and~\ref{lower}.

First, in Section~\ref{upper}, we prove an upper large deviations inequality or,
more precisely, the following exponential bound for the
probability that the chemical distance between two points $x$ and $y$ is
abnormally large:
\begin{theorem}
\label{devsup}
For every $p>p_c(d)$ and every $\varepsilon>0$,
we have
\[
\mathop{\lim\sup}_{\|x\|_1\to +\infty}\frac{\ln \P_p
 ( 0\communique x,   D(0,x)\ge (1+\varepsilon) \mu(x)  )}{\|x\|_1}<0.
 \]
\end{theorem}

The proof of this result strongly relies, through an appropriate renormalization
argument,  on the fact that, when $p$ is sufficiently close
to one, the chemical distance looks like the usual distance
$\|\cdot\|_1$.
\begin{theorem}
\label{grosseperco}
For each  $\alpha>0$, there exists $p'(\alpha)\in(p_c(d),1)$ such that for
every
 $p\in (p'(\alpha),1]$, the Bernoulli percolation with parameter $p$
 satisfies:
\[
\mathop{\lim\sup}_{\Vert x\Vert_1\to +\infty} \frac{\ln \P_p(0\communique
  x,   D(0,x)\ge (1+\alpha)\Vert x\Vert_1)}{\Vert x\Vert_1}<0.
  \]
\end{theorem}

We also obtain, as a corollary of this result, the continuity in $p=1$ of
the map $p\mapsto \mu_p$, where $\mu_p$ denotes the norm associated to the
chemical distance in the Bernoulli percolation with parameter $p$:
\begin{coro}
\label{conti}
$\lim_{p \to 1} \sup_{\|x\|_1\le 1}
 | \mu_p(x) - \|x\|_1  | =0$.
\end{coro}

In Section~\ref{lower}, we prove a lower large deviations inequality or,
more precisely, the following exponential bound for the
probability that the chemical distance between two points $x$ and $y$ is
abnormally small.
\begin{theorem}
\label{devinf}
For every $p>p_c(d)$ and every $\varepsilon>0$, we have
\[
\mathop{\lim\sup}_{\|x\|_1\to +\infty}\frac{\ln \P_p
 ( 0\communique x,   D(0,x)\le (1-\varepsilon) \mu(x)  )}{\|x\|_1}<0.
 \]
\end{theorem}

In its main lines, the proof follows the strategy used by Grimmett and
Kesten \cite{gk84} to prove an exponential bound for an analogous
quantity concerning first-passage percolation along the first coordinate axis.
However, two types of extra difficulties arise in our context: we want to
obtain an
exponential bound in every direction, not only along the first-coordinate
axis, and moreover we want this bound to be uniform
with respect to this direction. Thus, we first study in Lemma~\ref{tempsplan}
the minimal number of open edges needed to join the origin to hyperplanes
with a given direction. Then in Lemma~\ref{lesboites}
we study the minimal number of open edges needed to cross a box oriented
along the same direction. All estimates are done uniformly in the direction,
and, to conclude the proof of Theorem~\ref{devinf}, we use a renormalization argument.

Let us discuss briefly the speed---in $\|x\|_1$---that appears in the previous
large deviations inequalities.
Let us first look at the lower large deviations.  Choose an $x \in \Zd$, and then,
by the classical FKG inequalities, we obtain
\begin{eqnarray*}
&& \P_p  \bigl( D\bigl(0,(m+n)x\bigr) \le (1-\varepsilon) (m+n)\mu(x)  \bigr) \\
&& \qquad\ge  \P_p  \bigl( D(0,mx) \le (1-\varepsilon) m\mu(x),
D\bigl(mx,(m+n)x\bigr) \le (1-\varepsilon) n\mu(x) \bigr) \\
&&\qquad \ge  \P_p  \bigl( D(0,mx) \le (1-\varepsilon) m\mu(x)  \bigr)
\P_p  \bigl( D(0,nx) \le (1-\varepsilon) n\mu(x)  \bigr).
\end{eqnarray*}
Thus, the limit $\frac1{n\|x\|_1} \ln \P_p  ( D(0,nx)
\le (1-\varepsilon) n\mu(x)  )$ exists as $n$ goes to infinity,
and is strictly negative by Theorem~\ref{devinf}. Now, two distinct cases can occur:
        \begin{itemize}
        \item[$\bullet$] Either $\mu(x)=\|x\|_1$. This corresponds to the
existence of            a flat face in the asymptotic shape and occurs for
some $x$ as soon as $p>\overrightarrow{p_c}(d)$, critical probability for oriented
percolation on $\Zd$ (see \cite{gm04}).
In this case,
because of the inequality $D(0,x)\ge\|x\|_1$, the asymptotic  speed in the direction of $x$ is as fast as it is permitted by the geometry of the lattice: thus, we have
\[
\frac1{n} \ln \P_p  \bigl( D(0,nx) \le (1-\varepsilon) n\mu(x)  \bigr)=-\infty.
\]
        \item[$\bullet$] Or $\mu(x)>\|x\|_1$. Then for any $\varepsilon$ small enough
to have  $(1-\varepsilon)\mu(x)>\|x\|_1$,  we can force a deterministic path with
exactly $n\|x\|_1$ edges with ends $0$ and $nx$ to be open, which implies
that the chemical distance between $0$ and $nx$ is less than $(1-\varepsilon) n\mu(x)$:
\begin{eqnarray*}
&& \P_p  \bigl( D(0,nx) \le (1-\varepsilon) n\mu(x)  \bigr) \ge p^{n\|x\|_1}, \\
& \mbox{that is, } & \frac1{n \|x\|_1} \ln \P_p  \bigl( D(0,nx)
\le (1-\varepsilon) n\mu(x)  \bigr)
\ge \ln p >-\infty.
\end{eqnarray*}
        \end{itemize}
Thus the exponential rate in Theorem~\ref{devinf} is optimal.

Turning to the upper large deviations, we can once again force a deterministic path with
exactly $ \lfloor n(1+\varepsilon)\|x\|_1 \rfloor +1$ edges with ends
$0$ and $nx$ to be the only open path between $0$ and $nx$, which implies
that the chemical distance between $0$ and $nx$ is larger than
$(1+\varepsilon) n\mu(x)$:
\begin{eqnarray*}
\P_p  \bigl(+\infty> D(0,nx) \ge (1+\varepsilon) n\mu(x)  \bigr)
\ge  \bigl( p(1-p)^{2d}  \bigr)^{\lfloor n(1+\varepsilon)\|x\|_1 \rfloor +1},
\end{eqnarray*}
whence
\begin{eqnarray*}
\mathop{\lim\inf}_{n\to +\infty}
\frac{\ln \P_p  ( +\infty >D(0,nx) \ge (1+\varepsilon) n\mu(x)  )}{n \|x\|_1}
\ge (1+\varepsilon)\ln  \bigl( p(1-p)^{2d}  \bigr).
\end{eqnarray*}
Once again, the exponential rate  in Theorem~\ref{devsup} is optimal.

This phenomenon is quite different from what is expected of large deviations in
first-passage percolation with bounded passage times. Indeed, in the context of first-passage percolation with bounded passage times,
building a bad configuration that forces the travel time between $0$ and $nx$ to
be too large should cost more that~$c^{n \|x\|_1}$.
We expect then a speed in $\|x\|_1^d$; see \cite{k86} and also \cite{cz03}.
On the other hand, building a configuration that allows the
travel time between $0$ and $nx$ to
be too small should typically still need a cost of order $c^{n \|x\|_1}$,
as it is sufficient to build one ``too good'' path. Thus the speeds for upper
large deviations and lower large deviations in classical first-passage
percolation could be different.

Finally, in Section~\ref{forme}, thanks to the uniformity with respect to
the direction provided by Theorem~\ref{ap-enrhume}, we will also prove
a large deviation inequality
 for the asymptotic shape of the set $B_t$ of points that are at a distance
less or equal to $t$ from the origin:
\[
B_t=\{x \in \Zd\dvtx   0 \communique x,   D(0,x) \leq t \}.
\]
Since $p> p_c$, we can condition the probability measure on the event that
the origin~$0$
is in an infinite cluster, which has positive probability:
\[
\Pcond_p(A)=\frac{\P_p(A\cap \{0\in C_{\infty}\})}{\P_p(0\in
  C_{\infty})}.
  \]
In order to study the convergence of the random set $B_t/t$, we also introduce
the Hausdorff distance between two non empty compact subsets  of $\Rd$:
\begin{longlist}[2.]
\item[1.] For $x\in\Rd$ and $r\ge 0$, $\mathcal{B}_{\mu}(x,r)=
\{y\in\Rd\dvtx   \mu(x-y)\le r\}.$
\item[2.] The Hausdorff distance between two nonempty compact subsets  $\mathcal{K}_1$ and
$\mathcal{K}_2$ of $\Rd$
is defined by
\[
\mathcal{D}(\mathcal{K}_1,\mathcal{K}_2)=\inf\{r\ge 0\dvtx   \mathcal{K}_1\subset \mathcal{K}_2+\mathcal{B}_{\mu}(0,r)\mbox{ and }
\mathcal{K}_2\subset \mathcal{K}_1+\mathcal{B}_{\mu}(0,r)\}.\]
\end{longlist}
Note that the equivalence
of norms on $\Rd$ ensures that the topology induced by this Hausdorff distance
does not depend on the choice of the norm $\mu$. We can now state the random set
version of Theorem~\ref{ap-enrhume}.
\begin{theorem}
\label{dev-AS}
For every $p>p_c(d)$, for every $\epsilon>0$, there exist two strictly
positive constants $A$ and $B$ such that
\[
\forall t>0 \qquad \Pcond_p  \biggl(
\mathcal{D} \biggl(\frac{B_t}t,\mathcal{B}_{\mu}(0,1) \biggr) \ge \epsilon
 \biggr) \le A e^{-Bt}.
\]
\end{theorem}

This result improves the following asymptotic shape result that was proved by the
authors in \cite{gm04}: for every $p>p_c(d)$,
\[
 \lim_{t\to +\infty}\mathcal{D} \biggl(\frac{B_t}t,\mathcal{B}_{\mu}(0,1) \biggr)=
0,\qquad\Pcond_p\as
\]

Let us begin now with the main notations and a reminder of some common
useful results in supercritical percolation theory. We also include in the
following section a technical lemma to build bases of $\Rd$ that are adapted
to the proof of directional
estimates.

\section{Notations and preliminary results}
\label{notation}

\subsection{Norms\textup{,} balls and spheres}
On the space $\Rd$, consider the canonical basis $(e_1,\dots, e_d)$. For every $x \in \Rd$,
define the three classical following norms:
\[
\|x\|_1 =  \sum_{m=1}^d |x_i|,\qquad
\|x\|_2 =   \Biggl(\sum_{m=1}^d |x_i|^2 \Biggr)^{1/2},\qquad
\|x\|_\infty =  \max_{1 \le m \le d} |x_i|.
\]
For $i \in \{1,2,\infty\}$, $x \in \Rd$ and $r>0$, we define the following balls
in $\Zd$:
\[
\mathcal{B}_i(x,r)=\{y \in \Zd\dvtx  \|y-x\|_i \le r \}
\quad\mbox{and}\quad
\mathcal S_i=\{x \in \Rd\dvtx   \|x\|_i =1\}.
\]
Recall that
\[
\|x\|_\infty \le \|x\|_2 \le \|x\|_1 \le \sqrt d \|x\|_2.
\]
We also consider the norm $\mu$ given by Proposition~\ref{puiskilfo}. We recall the reader that, for $x \in \Rd$ and $r>0$, we chose
to consider, for the norm $\mu$, balls in $\Rd$ rather than in $\Zd$:
\[
\mathcal{B}_\mu(x,r)=\{y \in \Rd\dvtx   \mu(y-x)\le r \}.
\]

We also introduce $\mu_{\mathrm{inf}}= \mathop{\inf}_{y \in \mathcal S_1} \mu(y)$,
which is strictly positive. As $\mu$ is invariant under the symmetries of the grid,
we get the inequality
\[
\mu_{\mathrm{inf}}\|x\|_1 \le \mu(x)\le \mu(e_1)\|x\|_1.
\]

\subsection{Exponential inequalities}
Let us rewrite the result of Antal and Pisztora~\cite{ap96} in an appropriate form
to further computations: there exist three strictly positive constants $A_1$,
$B_1$ and $\rho$, depending only on the dimension $d$ and on the percolation
parameter $p>p_c(d)$, such that
\begin{equation}
\label{antal_pisztora}
\forall x\in\Zd \qquad
\P_p  \bigl( 0\communique x,   D(0,x)\ge \rho\|x\|_1  \bigr)
\le A_1e^{-B_1\|x\|_1}.
\end{equation}

We also recall here some classical results concerning the geometry of the clusters
in supercritical percolation.
Thanks to \cite{ccgks89}, we
can control the radius of finite clusters: there exist two strictly positive
constants $A_2$ and $B_2$ such that
\begin{equation}
\label{amasfini}
\forall r>0 \qquad
 \P_p  \bigl( |C(0)|<+\infty,   0\communique\partial \mathcal B_1(0,r)  \bigr)
\le A_2e^{-B_2 r}.
\end{equation}

We can also control the size of the holes in the infinite cluster: there exist
two strictly positive constants $A_3$ and $B_3$ such that
\begin{equation}
\label{amasinfini}
\forall r>0 \qquad
 \P_p  \bigl( C_\infty \cap \mathcal B_1(0,r) =\varnothing  \bigr)
\le A_3e^{-B_3 r}.
\end{equation}
When $d=2$, this result follows from the large deviation estimates by Durrett
and Schonmann~\cite{ds88}. Their methods can easily be transposed when $d\ge 3$.
Nevertheless, when $d\ge 3$, the easiest way to obtain it seems to use \cite{gm90} slab's result.

Note that in (\ref{amasfini}) and in (\ref{amasinfini}), the choice of
the norm $\|\cdot\|_1$ is, of course, irrelevant thanks to the norm equivalence.

\subsection{Stochastic comparison}
First, there is a natural partial order $\preceq$ on $\Omega=\{0,1\}^{\Ed}$: for
$\omega$ and $\omega'$ in $\Omega$, one says that $\omega\preceq\omega'$ holds
if and only if $\omega_e\le \omega'_e$ for each $e\in\Ed$.
Consequently, we say that a function $f\dvtx\omega\to\R$ is nondecreasing
if $f(\omega)\le f(\omega')$ as soon as $\omega\preceq\omega'$.
An event $A$ is said to be nondecreasing if its indicator function $\1_A$ is
nondecreasing.

Let us now recall the concept of stochastic domination: we say that a probability
measure $\mu$ dominates a probability measure $\nu$ if
\[
{\int f \, d\nu}\le{\int f \, d\mu}
\]
holds as soon as $f$ in an nondecreasing function.
We also write $\nu\preceq\mu$.

In the following, it will often be useful  to compare locally dependent fields with
products of Bernoulli probability measures: remember that a family $\{Y_x,x\in\Zd\}$
of random variables is said to be locally dependent if there exists $k$ such that,
for every $a \in \Zd$, $Y_a$ is independent of $\{Y_x\dvtx \|x-a\|_2\ge k\}$.

\begin{prop}[(\cite{lss97})]
\label{comparaison}
Let $d,k$ be positive integers. There exists a nondecreasing function $\pi\dvtx[0,1]\to[0,1]$
satisfying $\lim_{\delta\to 1}\pi(\delta)=1$ such that the following holds:
 if $Y=\{Y_x,x\in\Zd\}$ is a locally dependent family of random variables satisfying
 \[
 \forall x\in\Zd \qquad P(Y_x=1)\ge\delta,
 \]
then
$P_Y\succeq \Ber(\pi(\delta))^{\otimes\Zd}$.
\end{prop}

This is in fact a particular case, but sufficient for our purposes, of a more
general result given in \cite{lss97}.

\subsection{Some consequences of the symmetry properties}
Let us introduce some notation: we denote by $\Sd$ the symmetric group on $\{1,\dots,d\}$.
For each $x=(x_1,\dots,x_d)\in\Rd$, $\sigma \in \Sd$ and $\epsilon\in\{+1,-1\}^d$, we define
\[
\Psi_{\sigma,\epsilon}(x)=\sum_{i=1}^d\epsilon(i)x_{\sigma(i)}e_i.
\]
Then $\mathcal O(\Zd)=\{\Psi_{\sigma,\epsilon}\dvtx    \sigma \in \Sd, \, \epsilon\in\{+1,-1\}^d\}$ is the
group of orthogonal transformations that preserve the grid $\Zd$. Consequently,
its elements also preserve the norm $\mu$.

When studying the chemical distance in a given direction $x$, we want to find a
basis of $\Rd$ adapted to the studied direction, that is, made of images of $x$ by
elements of $\mathcal O(\Zd)$. The next technical lemma gives the existence
of such a basis,
and an extra uniformity property in the direction $y$:

\begin{lemme}
\label{controle-direction}
There exists a constant $C_d>0$ such that, for each $x\in\Rd$, there exists a
family $(g_{1,x},g_{2,x},\dots, g_{d,x})\in(\mathcal{O}(\Zd))^d$ with
$g_{1,x}=\mathrm{Id}_{\Rd}$ and
 such that the
linear map $L_x :\Rd\to\Rd$ defined by
\[
\forall i\in\{1,\dots,d\} \qquad   L_x (e_i)=g_{i,x}(x)
\]
satisfies
\begin{equation}
\label{controle-L}
\forall y\in\Rd \qquad C_d\|y\|_1\|x\|_1\le\|L_x (y)\|_1\le\|y\|_1\|x\|_1.
\end{equation}
If moreover, for each  $n\in \mathcal{S}_2$, we set
$(n_1,n_2, \dots ,n_d)=(n, g_{2,n}( n), \dots,  g_{d,n}(n))$, then we
have
\begin{equation}
\label{controle-normales}
\forall y \in \Rd\qquad \frac{C_d}{d^{3/2}} \|y\|_2 \le
 \Biggl(
    \sum_{m=1}^d \langle y,n_m\rangle ^2  \Biggr) ^{1/2}
\le \sqrt d \|y\|_2.
\end{equation}
\end{lemme}

Note that in dimension two, this construction is much simpler: if $R$ denotes the rotation with angle $\pi/2$, we can set
\[
L_x(e_1)=x \quad\mbox{and}\quad L_x(e_2)=R(x).
\]
However, this is more intricate in higher dimension. For instance, in dimension three, if $x=(1,1,1)$, none of the images of $x$ by $\mathcal O(\Zd)$ is orthogonal to $x$. In particular, even if $\|x\|_2=1$, $L_x$ may not be in $\mathcal O (\Rd)$.

\begin{pf*}{Proof of Lemma \protect\ref{controle-direction}}
Choose $x \in \Zd$.
Then, for every $(g_1,\dots,g_d) \in  ( \mathcal{O}  ( \Zd  ) )^d $,
denote by  $A^x_{g_1,\dots,g_d}$ the only linear map which satisfies
\[
\forall i\in\{1,\dots,d\}\qquad A^x_{g_1,\dots,g_d}(e_i)=g_i(x).
\]
Let $y=\sum\limits_{i=1}^{d} y_ie_i \in\Rd$: by linearity, we have
\begin{eqnarray*}
\|A^x_{g_1,\dots,g_d} (y)\|_1
& = &   \Biggl\|\sum_{i=1}^d y_i A^x_{g_1,\dots,g_d} (e_i) \Biggr\|_1
=  \Biggl\|\sum_{i=1}^d y_i g_i(x) \Biggr\|_1\\
& \le & \sum_{i=1}^d |y_i| \times \|A^x_{g_1,\dots,g_d} (e_i)\|_1
=  \sum_{i=1}^d |y_i| \times \|x\|_1=\|y\|_1\|x\|_1.
\end{eqnarray*}
Now define, for each $x \in \Rd$,
\begin{eqnarray*}
b(x) & = & \max_{(g_2,\dots,g_d)\in \mathcal{O}(\Zd)^{d-1}}
\inf_{y\in\mathcal{S}_1}\|A^x_{\mathrm{Id},g_2,\dots,g_d} (y)\|_1,
\end{eqnarray*}
and define $L_x$ to be an application $A^x_{g_1,\dots,g_d}$ which realizes
the maximum in the definition of $b(x)$.
Let us set
\[
C_d  =  \inf_{x\in\mathcal{S}_1} b(x).
\]
It is easy to see that, for every $x \in \Rd$, $L_x$ satisfies
equation~(\ref{controle-L}) and it only remains to prove that $C_d>0$.

Clearly, $x\mapsto b(x)$ is a continuous map.
So, since $\mathcal S_1$
 is a compact set, it is sufficient to prove
that $b(x)\ne 0$ for any $x \in \mathcal S_1$.
Let then $x\ne 0$: there exists $i_0$ such that $x_{i_0}\ne 0$.
Consider $i\in\{1,\dots, d\}$; we can find $\sigma\in\Sd$ with $\sigma (i)=i_0$.
Now let $h\in\{-1,+1\}^d$ with $h(i)=-1$ and $h(j)=1$ for $i\ne j$: then
one has $\Psi_{\sigma,(1,\dots,1)}(x)- \Psi_{\sigma,h}(x)=2x_{i_0}e_i$.
It follows that the vector space generated by
$\{g(x)\dvtx   g\in \mathcal{O}(\Zd)\}$ is equal to $\Rd$.
Then, since $x\ne 0$, one can find a family
$(g_2,\dots,g_d)\in (\mathcal{O}(\Zd))^{d-1}$ such that
$(x,g_2(x),\dots,g_d(x))$ is a basis of $\Rd$.
Thus, $A^x_{g_1,\dots,g_d}$ is a linear invertible map.
This implies that $\inf_{\|y\|_1=1}\|A^x_{g_1,\dots,g_d} (y)\|_1>0$, and
hence that $b(x)>0$.

Let us prove inequality~(\ref{controle-normales}). If we define
\[
B(x)=\sum_{m=1}^d \langle  n_m, x\rangle e_m,
\]
then we have $\sum_{m=1}^d \langle x,n_m \rangle^2=\|B(x)\|_2$, and moreover,
$\langle B (e_j),e_i\rangle=\langle n_i,e_j\rangle$
$=\langle L_{ n} (e_i),e_j\rangle=\langle L_{n}^{*} (e_j),e_i\rangle$,
which is equivalent to say that $B=L_{n}^{*}$.
Equation~(\ref{controle-L}) and the equivalence of norms imply then that
\begin{equation}
\label{from-below}
\forall x\in\Rd\qquad \| L_{n} (x)\|_2\ge \frac{C_d}{d^{3/2}}\|x\|_2.
\end{equation}
Let us denote by $  |\hspace{-1.51pt}\|A |\hspace{-1.51pt}\|_2=\sup_{x\in\mathcal{S}_2} \|Ax\|_2$.
We have:
\[
\| L_{n}^{*} x\|_2 \ge \frac1{ |\hspace{-1.51pt}\|(L_{ n}^{*})^{-1} |\hspace{-1.51pt}\|}_2\|x\|_2.
\]
It is clear from (\ref{from-below}) that
$ |\hspace{-1.51pt}\| L_{n}^{-1}|\hspace{-1.51pt}\|_{2} \le  \frac{d^{3/2}}{C_d}$.
Applying to $A= L_{n}^{-1}$ the classical identity
\[
 |\hspace{-1.51pt}\| A |\hspace{-1.51pt}\|_2 =\sup_{x\in\mathcal{S}_2} \|Ax\|_2
=\sup_{x\in\mathcal{S}_2,y\in\mathcal{S}_2} \langle Ax,y\rangle
=\sup_{x\in\mathcal{S}_2} \|A^{*}x\|_2=  |\hspace{-1.51pt}\|A^*|\hspace{-1.51pt}\|_2,
\]
it follows that
\[
\forall x\in\Rd\qquad \Biggl (\sum_{m=1}^d \langle x,n_m \rangle^2 \Biggr)^{1/2}
=\| L_{ n}^{*} x\|_2\ge \frac{C_d}{d^{3/2}}\|x\|_2,
\]
which concludes the proof of the left-hand side. The right-hand side is obvious.
\end{pf*}

\section[Upper large deviations: Proof of Theorem~1.3]{Upper large deviations: Proof of Theorem~\protect\ref{devsup}}
\label{upper}
The aim of this section is to prove the upper large deviations estimate,
Theorem~\ref{devsup}, for the chemical distance. First, we prove the exponential
inequality for $p$ close to $1$ given by Theorem~\ref{grosseperco}, then we deduce Corollary~\ref{conti} and finally, via a renormalization argument,
we prove the large deviations result for every $p>p_c$.

\subsection[Chemical distance for $p$ close to $1$\textup{:} proof of Theorem~1.4]{Chemical
distance for $p$ close to $1$\textup{:} proof of
Theorem~\textup{\protect\ref{grosseperco}}} %
For this proof, we also consider the $*$-topology on $\Zd$: two points $x,y \in \Zd$
are $*$-\textit{neighbors} if and
only if $\|x-y\|_\infty=1$. A $*$-\textit{path} is a sequence $(x_1,\dots,x_n)$ such
that for
every $i \in \{1,\dots,n-1\}$, $x_i$ and $x_{i+1}$ are $*$-neighbors. A
set $E$ is
$*$-connected if between any two of its
vertices, there exists a $*$-path using only vertices in $E$.

Given a configuration $\omega$, say that a point $x\in\Zd$
is \textit{wired} if each bond $e=(s,t)$ with $\|s-x\|_{\infty}\le 1$ and
$\|t-x\|_{\infty}\le 1$ satisfies $\omega_e=1$. Otherwise, say that $x$ is
\textit{unwired}. The wired points should be considered as the good guys,
whereas the unwired points are the bad ones.
Let $Y_x=\1_{\{x\mathrm{\, is\, unwired}\}}$; then $x$ is wired if and only if
$Y_x=0$.

Let us define $V(x)(\omega)$ to be the set of points $y\in\Zd$
such that there exists a $*$-path of unwired vertices from $x$ to $y$, which
means that there
exist $n\ge 0$ and $x=x_0,x_1,\dots,x_n=y$,
with $Y_{x_i}=1$ for each $i\in\{0,\dots,n\}$ and
$\|x_i-x_{i+1}\|_{\infty}=1$ for each $i\in\{0,\dots, n-1\}$.
Note that $V(x)=\varnothing$ as soon as $x$ is wired.
By definition,~$V(x)$ is always a $*$-connected set.
For $x\in\Zd$, we define 
\[
V_1(x)=V(x)+\{-1,0,1\}^d\quad\mbox{and}\quad V_2(x)=V_1(x)+\{-1,0,1\}^d.
\]

Let us show that when $p$ is large enough, $V(x)$ is almost surely a finite set.
For each $p\in[0,1]$, the field $(Y_x)_{x\in\Zd}$ is a locally dependent
$\{0,1\}$ valued stationary field, with $\lim_{p\to 1}\P_p(Y_0=1)=1$.
It follows from Proposition~\ref{comparaison} that there exists~$r_1(p)$ with
\[
(\P_p)_Y\succeq
\Ber(r_1(p))^{\otimes\Zd}
\]
and $\lim_{p\to 1} r_1(p)=1$. We can thus find $p'_1 \in
(p_c(d),1)$ such that for every \mbox{$p>p'_1$},
\begin{equation}
\label{p'1}
(2d-1)\bigl(1-r_1(p)\bigr)<1.
\end{equation}
By a classical counting argument, this ensures that $V(x)$ is $\P_p$ almost
surely a finite set. Suppose for the sequel that $p>p'_1$. Under this
assumption, we have the following result:

\begin{lemme}
Let $x \in \Zd$. Suppose that $s,t\in V_1(x)$ with $s\communique t$. Then, there
exists an open path  from $s$ to $t$ which only uses vertices in
$V_2(x)$.
\end{lemme}

\begin{pf}
Since $V(x)$ is bounded, $V(x)^c$ has only finitely many $*$-connected components
and exactly one of them is of infinite size.
Of course, a path from  $s$ to $t$ can meet one or more of these sets.
We will prove that for every connected component $K$ of $V(x)^c$ and every
 open path $\beta$ from $s$ to $t$, the path $\beta$ can be modified to get an open path
from $s$ to $t$ which never enters $K\backslash V_2(x)$.

Suppose that $\beta=(s=x_0,x_1,\dots,x_n=t)$ and define
$i=\min\{k\ge 0\dvtx   x_k\in K\}$ and $j=\max\{k\ge 0\dvtx    x_k\in K\}$. Clearly $i>0$ and
$j<n$. Obviously, $\{x_i,x_j\}\subset\frinstar(K)$, where $\frinstar(K)$ is
the set of points $x$ in $K$ such that there exists $y \in \Zd \backslash
K$ with $\|x-y\|_\infty=1$.
By part (ii) of Lemma 2.1 in \cite{dp96},
the set $\frinstar(K)$ is $*$-connected. Note that by definition of $V(x)$, every point of
$\frinstar(K)$ is wired.
But if $a$ and $b$ are $*$-neighbors, then $b \in a+\{-1,0,1\}^d$. Since
$a$ is wired, there exists an open path from $a$ to $b$ using only edges in
$a+\{-1,0,1\}^d$. Then, there exists a open path from $x_i$ to $x_j$ which only
uses points in  $\frinstar(K)+\{-1,0,1\}^d\subset V_1(x)+\{-1,0,1\}^d=V_2(x)$.
\end{pf}

We can now come back to the proof of the theorem. Let $\alpha>0$.
Choose $x$ in $\Zd$ and let $\gamma=(0=x_0,x_1,\dots,x_n=x)$ be a fixed path from
$0$ to $x$ with the minimal possible number of edges $n=\|x\|_1$. We let
\[
V=\bigcup_{y\in\gamma} V(y).
\]
Now suppose that there exists an open path from $0$ to $x$. Let us prove that
under this condition, we can find an open path from $0$ to $x$
which only uses points in $\gamma\cup(V+\{-2,-1,0,1,2\}^d)$.

Let $i$ be the greatest
integer in $\{0,\dots,n\}$ such that  there exists an open path from $0$ to
$x_i$ which only uses points in
$\gamma\cup(V+\{-2,-1,0,1,2\}^d)$. Note that since $0\communique x_i$
and $0\communique x$, we have $x_i\communique x$.  We want to prove that $i=n$.

Suppose by contradiction that $i<n$.
The maximality of $i$ implies that $x_i$ can not be wired. So $V(x_i)\supset\{x_i\}$,
therefore it
is not empty, which allows to define $j=\max\{k\in\{i+1,\dots,n\}\dvtx   x_k\in V(x_i)\}$.
\begin{itemize}
\item[$\bullet$] If $j=n$, then $x_i$ and $x$ belong to $V(x_i)$. Since
$x_i\communique x$, it follows from the previous lemma that
there exists an open  path from $x_i$ to $x$ which only uses vertices in
$V_2(x_i)$. Joint with the part of the path $\gamma$ from $0$ to $x_i$, this
gives  an open path from $0$ to $x$ which only uses points that are in
$\gamma\cup(V+\{-2,-1,0,1,2\}^d)$, which is a contradiction.
\item[$\bullet$]
If $j<n$, let $K$ be the connected component of $V(x_i)^c$ which contains $x$.
On one side, $x_{i}\notin K$, so there exists $l\in\{i+1,\dots ,n\}$ such that
$x_l\in\frinstar(K)$.
On the other side, $x_{i} \notin K$ and $x\communique x_{i-1}$ by
definition of $i$, thus there exists $z\in\frinstar(K)$ such that  $x\communique z$.
Since points in $\frinstar(K)$ are linked, $x_l\communique z$, but $x
\communique z$ and $x_{i} \communique x$ so finally
$x_{i}\communique x_l$.
Since $\{x_i,x_l\} \subset V_1(x_i)$, by the previous lemma we see
that there exists an open path from $x_i$ to $x_l$ using only points of $V_2(x_i)$,
which contradicts again the maximality of $i$.
\end{itemize}
Thus under the assumption that
$0$ and $x$ belong to the same cluster, we have constructed an open path from $0$
to $x$
which only uses points in $\gamma\cup(V+\{-2,-1,0,1,2\}^d)$, and thus is not too
far away from the deterministic path $\gamma$.

Define, for every $y \in \Zd$, the event
\[
F_y=\bigcup_{z\dvtx   \|z-y\|_\infty \le 2}
\{z\mbox{ is unwired}\}
\]
and set
 $Z_y=\1_{F_y}$.
Since  $(Z_y)_{y\in\Zd}$ is locally dependent with
\mbox{$\lim_{p\to 1}\P_p(Z_{y}=1)=0$}, it follows from Proposition~\ref{comparaison} that
there exists $r_2(p)$ with
\[
(\P_p)_Z\preceq \Ber(r_2(p))^{\otimes\Ed}
\]
and $\lim_{p\to 1} r_2(p)=0$.
Note that by definition of $(Z_y)_{y\in\Zd}$, the open path we built
only uses points $y$ that are in $\gamma$ or satisfy $Z_y=1$.
Moreover, if we suppose now that $D(0,x)\ge (1+\alpha)\Vert x\Vert_1$, the
length of
this path is also greater than  $(1+\alpha)\Vert x\Vert_1$.

The idea is now the following: if $0 \communique x$ and
$D(0,x) \ge (1+\alpha)\|x\|_1$, then by the previous construction, there
must exist an open path between $0$ and $x$ with length larger than
$(1+\alpha)\|x\|_1$ and that
contains only points in $\gamma\cup(V+\{-2,-1,0,1,2\}^d)$: this path must then
contain at least $\alpha\|x\|_1$ points such that $Z_y=1$, which is unlikely
when $p$ is large enough. Let us turn this crude argument into a rigorous proof
via a counting argument.

Let $\Gamma$ be the family of self-avoiding paths from $0$ to $x$.
Clearly, if $\beta \in \Gamma$,
\[
\P_p(\forall y\in\beta\backslash\gamma, \, Z_y=1)\le r^{|\beta|-|\gamma|},
\]
where $r=r_2(p)$ and $|\beta|$ denotes the number of edges in
$\beta$. Remember that $\lim_{p\to 1} r_2(p)=0$.
We can thus find $p'_2(\alpha) \in (p_c(d),1)$ such that $\forall
p>p'_2(\alpha)$, $r=r_2(p)$ satisfies
\begin{equation}
\label{p'2}
(2d-1)r<1\quad\mbox{and}\quad (2d-1)^{1+\alpha}r^{\alpha}<1.
\end{equation}
It follows that
\begin{eqnarray*}
&& \P_p\bigl(0\communique x,  D(0,x)\ge (1+\alpha)\Vert x\Vert_1\bigr) \\
&&\qquad \le  \sum_{\beta\in\Gamma\dvtx   |\beta|\ge (1+\alpha)|\gamma|} \P_p(\forall
i\in\beta\backslash\gamma,   Z_i=1)\\
&&\qquad \le   \sum_{\beta\in\Gamma\dvtx   |\beta|\ge (1+\alpha)|\gamma|} r^{|\beta|-|\gamma|}\\
&&\qquad \le  \sum_{n=(1+\alpha)|\gamma|}^{+\infty} (2d)(2d-1)^{n-1} r^{n-|\gamma|}\\
&&\qquad \le   \frac{2d r^{-|\gamma|}}{(2d-1)(1-(2d-1)r)}
 \bigl((2d-1)r \bigr)^{(1+\alpha)|\gamma|} \\
&&\qquad \le   \frac{2d }{(2d-1)(1-(2d-1)r)}  \bigl( (2d-1)^{1+\alpha}r^\alpha
   \bigr)^{|\gamma|}.
\end{eqnarray*}
As $|\gamma|=\|x\|_1$, taking $p'(\alpha)=\max\{p'_1(\alpha),p'_2(\alpha)\}$---quantities
respectively defined in (\ref{p'1}) and (\ref{p'2})---ends the proof
of the theorem.

\subsection[Continuity of $\mu_p$ at $p=1$\textup{:} Proof of Corollary~1.5]{Continuity of $\mu_p$ at $p=1$\textup{:} Proof of
Corollary~\textup{\protect\ref{conti}}}%
Let $\alpha>0$ and $x\in\Zd$. Using the Borel--Cantelli lemma and
Theorem~\ref{grosseperco}, we obtain for $p>p'(\alpha)$:
\[
\mathop{\lim\sup}_{n\to +\infty}\1_{\{nx\communique
  0\}}\frac{D(0,nx)}{n}\le(1+\alpha)\|x\|_1, \qquad  \Pcond_p\as
  \]
By the very definition of $\mu_p$, the left-hand side is equal to $\mu_p(x)$.
Using moreover the fact that $\mu_p(x)\ge \|x\|_1$, we have proved
\[
\forall\alpha>0, \, \forall p>p'(\alpha), \, \forall x\in\Zd \qquad \|x\|_1\le
\mu_p(x)\le(1+\alpha)\|x\|_1.
\]
By homogeneity and continuity of $\mu_p$ and $\|\cdot\|_1$, we obtain the same
property for $x\in\Q^d$, and next for $x\in\Rd$:
\[
\forall\alpha>0, \, \forall p>p'(\alpha),\, \forall x\in\Rd \qquad
\|x\|_1\le
\mu_p(x)\le(1+\alpha)\|x\|_1,
\]
which ends the proof.

\subsection[Upper large deviations\textup{:} proof of Theorem~1.3]{Upper large deviations\textup{:} proof of Theorem~\textup{\protect\ref{devsup}}}
We can now prove the upper large deviations result Theorem~\ref{devsup}
for the chemical distance for every \mbox{$p>p_c$}.

\textit{Step \textup{1}}. Choice of constants.

Let $p>p_c(d)$ and $\varepsilon>0$ be fixed.

As $\mu$ is a norm, it is bounded away from $0$ on the compact set
$\mathcal{S}_1$,
and we can choose $\eta>0$  small enough to have
\begin{equation}
\label{eta}
\hspace*{10mm}
\forall \hat x \in \mathcal{S}_1   \qquad \biggl( 1+\frac{3\eta}{2\rho}  \biggr)
(1+\eta)^2 \mu(\hat x )+2\eta<(1+\epsilon) \mu(\hat x) \quad \mbox{and}
\quad \eta<\frac8{9}\rho,
\end{equation}
where $\rho$ has been defined in~(\ref{antal_pisztora}).
Note also
\begin{equation}
\label{alpha}
\alpha=\frac{\eta}{2\rho}.
\end{equation}

From now on, let us denote by $M$ a fixed integer which is such  that
$M\ge\frac{d}{\eta}\max(\frac{\mu(e_1)}2,\frac{\rho}{3C_d})$, where $C_d$
is the constant given by Lemma~\ref{controle-direction}.

For every $\hat x \in \R^d \cap \mathcal{S}_1$,  there exists $\hat r \in
\frac{1}{M}\Z^d \cap \mathcal{S}_1$  such that $\|\hat x-\hat r\|_1 \le
\frac{d}{2M}$.  Then, by the previous choice of $M$, one has
\begin{equation}
\quad
\label{M}
\|\hat x-\hat r\|_1 \le  \frac{d}{2M}\le\frac{C_d\alpha }3
\quad \mbox{and} \quad
|\mu(\hat x)-\mu(\hat r)| \le \mu(e_1)\|\hat x-\hat r\|_1\le\eta.
\end{equation}
Intuitively, $\hat r$ is a rational direction that approaches the ``real''
direction $\hat x$. Note that the result in Theorem~\ref {devsup} is
uniform in the direction, and the proof of this uniformity will use the fact
that $\frac{1}{M}\Z^d \cap \mathcal{S}_1$ is a finite set.

\textit{Step }2. Renormalization.

For $x\in\Zd$ and $N\in\N$, let us define the following set around $Nx$:
\[
\mathcal{I}_{x}^N
=  \bigl\{ y \in \mathcal{B}_1\bigl(Nx,\sqrt N\bigr)\dvtx   y\communique \partial
  \mathcal{B}_1(Nx,N)  \bigr\},
  \]
where $\partial  \mathcal{B}_1(Nx,N)=\{y\in\Zd\dvtx   \|Nx-y\|_1=N\}$.
We define the related random variable $I_{x}^N$ by
\[
I_{x}^N=\cases{Nx, & \quad if  $\mathcal{I}_{x}^N=\varnothing$,\cr
\operatorname{inf} \mathcal{I}_{x}^N, & \quad otherwise,}
\]
where $\operatorname{inf} \mathcal{I}_{x}^N $ denotes the point in $\mathcal{I}_{x}^N$
which is the closest to $Nx$.  If there are several, use for instance the
lexicographic order on the coordinates  to choose a unique point.
Note that:
\begin{itemize}
\item[$\bullet$] the random variable $I_{x}^N$ only depends on the states
  of the edges  in $\mathcal{B}_1(Nx,N)$.
\item[$\bullet$] $\|Nx-I_{x}^N\|_1\le \sqrt N$.
\end{itemize}
Since $I_{x}^N$ is close to $Nx$, the chemical distance $D(I_0^N,I_{x}^N)$
should be of the same order as $N\mu(x)$.  This is rigorously proved in the
following lemma:

\begin{lemme}
\label{unps}
The following results hold $\P_p$ almost surely\textup{:}
\begin{itemize}
\item[$\bullet$] For each $x\in\Zd$, $I_{x}^N\communique\infty$ for large $N$.
\item[$\bullet$] The sequence $(I_0^N)_{N\ge 1}$ is convergent.
\item[$\bullet$] For each $x\in\Zd$,
\[
\frac{ D(I_0^N,I_{x}^N)}N\to\mu(x).
\]
\end{itemize}
\end{lemme}

\begin{pf}
$\bullet$ The first assertion easily follows from Borel--Cantelli arguments.
At first, it follows from the exponential decay of the radius of finite
clusters---see equation (\ref{amasfini})---that $\P_p$ almost surely,
$\mathcal{I}_{x}^N=B(Nx,\sqrt N)\cap  C_{\infty}$ for large $N$.
The fact that  $B(Nx,\sqrt N)\cap  C_{\infty}$ is $\P_p$ almost surely  nonempty for large $N$
is now a consequence of (\ref{amasinfini}).

$\bullet$ Let us denote by $H$ the smallest element of $C_{\infty}$,
that is, the point in $C_\infty$  which is the closest to $0$ and among these points
if there are several,  the one which is the smallest for the lexicographic
order on the coordinates.  For large $N$, we have $H\in B(0,\sqrt N)$, so
$(I_0^N)_{N\ge 1}$ converges to $H$.

$\bullet$ If $x=0$, there is nothing to prove. Suppose then that $x\ne 0$.
By the previous point, the sequence with general term $D(I_0^N,I_x^N)/N$ has the same asymptotic
behavior as the sequence with general term $D(H,I_x^N)/N$.
It was proved in Lemma~5.7 of \cite{gm04} that
for every $\varepsilon>0$,
\[
\P_p \left( \matrix{
\exists M>0 \mbox{ }  \forall y \in \Z^d \cr
 (\|y\|_1 \geq M\mbox{ and
  }y\communique 0) \Longrightarrow |D(0,y)-\mu(y)| \leq \varepsilon \|y\|_1
}
 \left|\vphantom{\matrix{
\exists M>0 \mbox{ }  \forall y \in \Z^d \cr
 (\|y\|_1 \geq M\mbox{ and
  }y\communique 0) \Longrightarrow |D(0,y)-\mu(y)| \leq \varepsilon \|y\|_1
}}\right.0
\communique \infty
 \right)=1.
 \]
  By taking $M=|C(0)|$ when $0 \not\communique \infty$, we can remove the conditioning and obtain
\[
\P_p  \pmatrix{
\exists M>0 \mbox{ }  \forall y \in \Z^d \cr
(\|y\|_1 \geq M\mbox{ and
  }y\communique 0) \Longrightarrow |D(0,y)-\mu(y)| \leq \varepsilon \|y\|_1
}=1.
 \]
 Using translation invariance, this implies
\[
\P_p  \pmatrix{
\forall x\in\Zd \mbox{ }  \exists M_x>0 \mbox{ }  \forall y \in \Z^d \cr
(\|y\|_1 \geq M_x\mbox{ and
  }y\communique x) \Longrightarrow |D(x,y)-\mu(y)| \leq \varepsilon \|y\|_1
}=1.
\]
This implies
\[
\lim_{N \to \infty} \frac{ D(I_0^N,I_x^N)-\mu(I_x^N)}{\|I_x^N\|_1} = 0, \qquad \P_p \as
\]
Since  $\|Nx-I_x^N\|_1\le \sqrt N$, one has
$|N \mu(x)-\mu(I_x^N)| \le \sqrt N\mu(e_1)$ and $\|I_x^N\|_1 \sim N \|x\|_1$; the desired result follows.
\end{pf}

\begin{figure}

\includegraphics{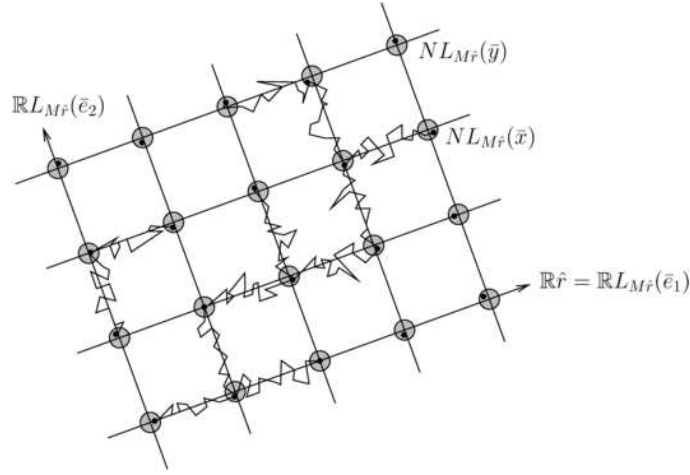}

\caption{The macroscopic grid $NL_{M\hat r} (\Zd)$ adapted to the
direction $\hat r$. The vertex $\bar x$ of the macroscopic grid has coordinates
$NL_{M\hat r} (\bar x)$ in the microscopic grid. It is surrounded by a
small ball with radius $\sqrt N$ which contains a black point, denoted by
$I^N_{L_{M\hat r}(\overline x)}$, which is  the  point of the infinite cluster
that is the closest to $NL_{M\hat r} (\bar x)$ in the microscopic grid.
The macroscopic edge between $\bar x$ and $\bar y$ is said to be open
if $D  ( I^N_{L_{M\hat
       r}(\overline x)},I^N_{L_{M\hat r}(\overline y)}  ) \le
  NM\mu(\hat r)(1+\eta)$, which happens with a probability going to $1$
  as $N$ goes to infinity. As the states of the macroscopic edges are locally dependent,
  by choosing $N$ large enough, the chemical distance in the macroscopic grid can be
  made very close to the $l^1$-distance.}
  \label{f1}
\end{figure}

We can now introduce a macroscopic  percolation: in order to study the chemical
distance in the direction $\hat r$, we are going to build a large grid,
with mesh~$NM$, whose axes are adapted to $
M\hat r$: the large grid is the image by $NL_{M\hat r}$ of the grid~$\Zd$,
where $L_{M\hat r}$ is given in Lemma~\ref{controle-direction} (see Figure \ref{f1}).
The macroscopic edge $\overline e=\{\overline x,\overline y\}$ has macroscopic
ends $\overline x$ and $\overline y$ that correspond in the microscopic graph
to the points $NL_{M\hat r}(\overline x)$ and $NL_{M\hat r}(\overline y)$;
for instance, the macroscopic vertex with coordinates $(1,0,\dots, 0)$ corresponds
to the vertex $NM\hat r$ in the microscopic lattice. By construction of $L_{M\hat r}$,
we expect the chemical distance between neighborhoods of the microscopic ends of any
macroscopic edge  to have a value of order $NM\mu(\hat r)$. If this event occurs,
we say that the corresponding macroscopic edge $\overline e$ is open, which should
happen with high probability.

\begin{lemme}
\label{grossearete}
For each $\hat r \in \frac{1}{M}\Z^d \cap \mathcal{S}_1$, for each positive integer $N$, we define
a field $(R^{N,\hat r}_{\overline e})_{\overline e\in\Ed}\dvtx$ if $\overline e=\{\overline x,\overline y\}$,
\[
R^{N,\hat r}_{\overline e}=\1_{G^{N,\hat r}_{\overline e}} \quad
\mbox{and} \quad G^{N,\hat r}_{\overline e}= \bigl\{ D  \bigl( I^N_{L_{M\hat
        r}(\overline x)},I^N_{L_{M\hat r}(\overline y)}  \bigr) \le
  NM\mu(\hat r)(1+\eta)  \bigr\}.
  \]
Then there exists a function $\overline{p}\dvtx\N\to [0,1]$,
independent of the choice of $\hat r \in \frac{1}{M}\Z^d \cap \mathcal{S}_1$  such that
\[
\lim_{N\to +\infty}\overline{p}(N)=1\quad \mbox{and}\quad
\P_{R^{N,{\hat r}}}\succeq\Ber(\overline{p}(N))^{\otimes\Ed},
\]
where $\P_{R^{N,{\hat r}}}$ denotes the law of the field
$(R^{N,\hat r}_{\overline e})_{e\in\Ed}$ on $\{0,1\}^{\Ed}$.
\end{lemme}

\begin{pf}
Our aim is  to apply once again the comparison result of Proposition~\ref{comparaison}.

Note that for any choice of $\hat r \in \frac{1}{M}\Z^d \cap
\mathcal{S}_1$, and any edge $\overline e=\{\overline x,\overline y\}$, the
macroscopic event $G^{N,{\hat r}}_{\overline e}$ only depends on states
of the microscopic edges in the ball
$\mathcal{B}_1 ( N L_{M\hat r}(\overline x), NM  (
    2+(1+\eta)\mu(e_1)  )  )$. Note also that, by
Lemma~\ref{controle-direction}, we have
\begin{eqnarray*}
&&\| \overline y- \overline x \|_1  > \frac2{C_d}  \bigl( 2+(1+\eta)\mu(e_1)  \bigr) \\
&&\qquad \Rightarrow
\| L_{M\hat r}(\overline y- \overline x) \|_1  > 2M \bigl( 2+(1+\eta)\mu(e_1)  \bigr) \\
 &&\qquad \Rightarrow
\| NL_{M\hat r }(\overline y)- NL_{M\hat r }(\overline x) \|_1  > 2NM \bigl( 2+(1+\eta)\mu(e_1)  \bigr).
\end{eqnarray*}
The field $(R^{N,{\hat r}}_{\overline  e})_{\overline e\in\Zd}$
is locally dependent, for some constant $k$ that does not depend on $N$, nor on the choice of $\hat r \in
\frac{1}{M}\Z^d \cap \mathcal{S}_1$, nor on $M$.
Since $(R^{N,{\hat r}}_{\overline  e})_{\overline e\in\Ed}$ is invariant
under translations
and symmetries of the grid $\Zd$, we only have to prove
that
\[
\lim_{N\to +\infty}\P_p  (G^{N,{\hat r}}_{\overline  e}
 )=1
 \]
 uniformly in $\hat r$ for the edge $\overline e=(0,e_1)$. But
the set $\frac{1}{M}\Z^d \cap \mathcal{S}_1$ is finite, so it is sufficient
to prove this limit for any $\hat r \in \frac{1}{M}\Z^d \cap \mathcal{S}_1$.
By applying Lemma~\ref{unps} to $x=L_{M\hat r}(e_1)$ for a given $\hat r \in \frac{1}{M}\Z^d \cap \mathcal{S}_1$ and
using the fact that almost sure convergence implies convergence in
probability, we end the proof of Lemma~\ref{grossearete}.
\end{pf}

Choose now $N$ large enough to be sure that $\overline{p}(N)$ given by
Lemma~\ref{grossearete} satisfies
\begin{equation}
\label{N}
\overline{p}(N)>p'  \biggl( \frac{1+3\alpha}{1+2\alpha}-1  \biggr)>p_c(d),
\end{equation}
where $p'(\cdot)$ is defined in Theorem~\ref{grosseperco} and $\alpha$ has been
defined in~(\ref{alpha}). For each $\hat r \in \frac1M \Zd \cap \mathcal{S}_1$,
we can construct a
macroscopic  percolation with mesh $NM$: we say that the edge $\overline e\in \Ed$
is open in the macroscopic  percolation associated to $\hat r$ if the event~$G_{\overline e}^{N, \hat r}$ occurs, and closed otherwise. This induces a
dependent percolation model on the macroscopic edges: all vertices in the macroscopic
lattice will be denoted by an overlined letter, the infinite cluster of
this macroscopic will be denoted $\overline{C}_\infty$, while the chemical distance
in this macroscopic lattice will still be denoted by $D$. The previous lemma compares
this locally dependent macroscopic percolation with i.i.d. Bernoulli percolation,
and the choice we made for $N$ allows us to use the result of Theorem~\ref{grosseperco}.

The strategy is now the following: for a large $x \in \Zd$, choose a $\hat r$
whose direction is close to the one of $x$  and build the macroscopic  percolation
associated to $\hat r$. Use Theorem~\ref{grosseperco} to find a macroscopic path
from a point not too far from $0$ to a point not too far from $x$, and whose length
is well controlled. Then come back to the initial microscopic percolation, and verify
that the existence of this macroscopic path implies, on the event $0 \communique x$,
the existence of an open microscopic path whose length is also well controlled.

\textit{Step }3. Construction of the macroscopic and microscopic paths.

From now on, we suppose, without loss of generality, that $x \in \Z^d$
satisfies
\begin{equation}
\label{x}
\|x\|_1 \geq \frac{8\sqrt N}\alpha.
\end{equation}
We emphasize that the constants $\alpha$ and $N$ have been defined in
(\ref{alpha}) and (\ref{N}) before
any choice of~$x$.

Then, we associate to $\hat x=x/\|x\|_1$ an approximate
$\hat r \in \frac{1}{M}\Z^d \cap \mathcal{S}_1$ satisfying equation~(\ref{M}).

We build the macroscopic  percolation associated to $\hat r$ and denote by
\begin{equation}
\label{xbar}
\bar x=\lfloor \|x\|_1/(NM) \rfloor  e_1
\end{equation}
 the vector in the coordinates of the macroscopic  grid that approximates $x$, where~$\lfloor t \rfloor$ denotes the integer part of the real number $t$.

Remember that, thanks to Lemma~\ref{controle-direction}, the application
$L_{M\hat r}$ maps $e_1$ to $M\hat r$  and satisfies
\[
 \forall i\in\{1,\dots,d\} \qquad  \mu(L_{M\hat r}(e_i))=M\mu(\hat{r})
 \quad\mbox{and}\quad  \|L_{M\hat r}(e_i)\|_1=M\]
 \begin{equation}
 \label{norme-de-L}
 \mbox{and }  \forall \hat r \in\frac{1}{M}\Z^d \cap \mathcal{S}_1\,\forall
 t\in\Rd\qquad C_d M\|t\|_1\le \|L_{M\hat r} (t)\|_1\le M\| t\|_1.
\end{equation}
For each $z\in\Zd$ and each $r>0$, we define the annulus
\[
\mathcal{A}(z,r)=\{y\in\Zd\dvtx   r/2\le \|y-z\|_1\le r\}
\]
and consider the following  ``good'' event of $\{0,1\}^{\Ed}$
in the macroscopic  percolation:
\[
G=
 \left\{
\matrix{
\exists \overline{a} \in \mathcal{A}(0,\alpha\|\overline{x}\|_1) \,\exists \overline{b}
\in \mathcal{A}(\overline{x},\alpha\|{\overline{x}}\|_1) \mbox{ such that }
\cr
\overline{a} \communiqueM \overline{b} \mbox{ and } D(\overline{a},\overline{b})
 \le (1+3\alpha)\|\overline{x}\|_1
}
 \right\}.
 \]
Note that for the complementary set of $G$, we have
\begin{eqnarray*}
G^c & \subset &
\{\mathcal{A}(0,\alpha\|\overline{x}\|_1)\cap
\overline{C}_{\infty}=\varnothing\}
\cup \{\mathcal{A}(\overline{x},\alpha\|\overline{x}\|_1)
\cap \overline{C}_{\infty}=\varnothing\}\\
& & \cup
\mathop{\bigcup_{\overline{a}\in
    \mathcal{B}_1(\overline{0},\alpha\|\overline{x}\|_1)}}_{\overline{b}\in
    \mathcal{B}_1(\overline{x},\alpha\|\overline{x}\|_1)}
\{\overline{a}\communiqueM \overline{b},    D(\overline{a},\overline{b})> (1+3\alpha)\|\overline{x}\|_1\} \\
& \subset & \{\mathcal{A}(0,\alpha\|\overline{x}\|_1)\cap \overline{C}_{\infty}=\varnothing\}  \cup \{\mathcal{A}(\overline{x},\alpha\|\overline{x}\|_1)\cap \overline{C}_{\infty}=\varnothing\}\\
& & \cup \mathop{\bigcup_{\overline{a}\in
    \mathcal{B}_1(\overline{0},\alpha\|\overline{x}\|_1)}}_{\overline{b}\in
    \mathcal{B}_1(\overline{x},\alpha\|\overline{x}\|_1)}
 \biggl\{\overline{a}\communiqueM \overline{b},   D(\overline{a},\overline{b})> \frac{1+3\alpha}{1+2\alpha} \|\overline{b}-\overline{a}\|_1 \biggr\}.
\end{eqnarray*}
As $G$ is an increasing event, we have by Lemma~\ref{grossearete}
\[
\P_{R^{N, \hat r}}(G^c)\le \Ber(\overline{p}(N))^{\otimes\Ed}(G^c).
\]
It follows that
\begin{eqnarray*}
\P_{R^{\epsilon,N}}(G^c) & \le & 2\P_{\overline{p}(N)}\bigl(\mathcal{B}_1(0,\alpha\|\overline{x}\|_1)\cap C_{\infty}=\varnothing\bigr)\\
&&+ \mathop{\sum_{\overline{a}\in
    \mathcal{B}_1(0,\alpha\|\overline{x}\|_1)}}_{\overline{b}\in
    \mathcal{B}_1(\overline{x},\alpha\|\overline{x}\|_1)}\P_{\overline{p}(N)}  \biggl( \overline{a}\communique \overline{b},   D(\overline{a},\overline{b})> \frac{1+3\alpha}{1+2\alpha}\|{\overline{a}-\overline{b}}\|_1
    \biggr).
\end{eqnarray*}
By the choice (\ref{N}) we made for $N$, the inequality
$\overline{p}(N)>p_c(d)$ is satisfied, so, by equation~(\ref{amasinfini}),
\[
 \P_{\overline{p}(N)}  \bigl(\mathcal{B}_1(0,\alpha\|\overline{x}\|_1)\cap C_{\infty}=\varnothing   \bigr)
\le A_3e^{-B_3\alpha\|\overline{x}\|_1}.
\]
Moreover, our choice of $N$ in (\ref{N}) was intended to apply Theorem~\ref{grosseperco}: thus, there exist two strictly positive constants $A$ and $B$ such that
\begin{eqnarray*}
&& \forall \overline{a}  \in \mathcal{B}_1(0,\alpha\|\overline{x}\|_1)
,\, \forall \overline{b} \in
\mathcal{B}_1(\overline{x},\alpha\|\overline{x}\|_1)  \\
&&\qquad
\P_{\overline{p}(N)}  \biggl( \overline{a}\communique \overline{b},   D(\overline{a},\overline{b})> \frac{1+3\alpha}{1+2\alpha}\|{\overline{a}-\overline{b}}\|_1   \biggr)
\le Ae^{-B\| \overline{b}-\overline{a}\|_1}.
\end{eqnarray*}
Thus we obtain
\begin{eqnarray*}
&& \P_{R^{\epsilon,N}}(G^c) \\
&&\qquad \le
2A_3e^{-B_3\alpha\|\overline{x}\|_1}+(C\alpha\|\overline{x}\|_1)^{2d}Ae^{-B (1-2\alpha)\|\overline{x}\|_1}\\
&&\qquad \le
2A_3e^{-B_3\alpha\|{x}\|_1/NM}+\bigl(C\alpha(\|x\|_1/NM+1)\bigr)^{2d}Ae^{-B (1-2\alpha)\|{x}\|_1/NM},
\end{eqnarray*}
where $C$ is a constant depending only on the dimension of the grid $\Zd$.

So, with a probability tending to $1$ exponentially fast with $\|x\|_1$, there exists a path in the macroscopic  percolation from a point in the set
$\mathcal{A}(0,\alpha\|\overline{x}\|_1)$ to  a point in the set
$\mathcal{A}(\overline{x},\alpha\|{\overline{x}}\|_1)$
which uses only edges $\overline e$ such that $G^{N, \hat r}_{\overline e}$ holds
and whose length is smaller or equal to $(1+3\alpha)\|\overline{x}\|_1$.
This implies the existence of a microscopic open path
from some 
microscopic vertex $S\in N L_{M\hat r}\mathcal{A}(0,\alpha\|\overline{x}\|_1)+\mathcal{B}_1(0, \sqrt N)$
to some 
microscopic vertex $T\in N L_{M\hat r}\mathcal{A}(\overline{x},\break \alpha\|\overline{x}\|_1)+\mathcal{B}_1(0, \sqrt N)$,
and whose length, by Lemma~\ref{grossearete}, is smaller than
$(1+3\alpha) \|\overline{x}\|_1 (1+\eta)NM \mu(\hat r) \le (1+3\alpha)(1+\eta)^2 \|{x}\|_1\mu(\hat x)$ by the choice of $\hat r$ in (\ref{M}) and the definition (\ref{xbar}) of $\overline{x}$.

\textit{Step }4. It remains now to link the ends $S$ and $T$ of this microscopic path to $0$ and~$x$ respectively, and to prove that with high probability, these links are very short.

Thanks to the definition of the annuli and to equations~(\ref{norme-de-L}) and~(\ref{x}), one has
\begin{eqnarray*}
\dfrac38{C_d\alpha\|x\|_1}&\le &\frac{C_d\alpha\|x\|_1}2-\sqrt N\le\|S\|_1\le\alpha\|x\|_1+\sqrt N\le  \frac98{\alpha\|x\|_1}, \\
\dfrac38{C_d\alpha\|x\|_1}&\le &\frac{C_d\alpha\|x\|_1}2-\sqrt N\le\|T-\|x\|_1\hat{r}\|_1\le\alpha\|x\|_1+\sqrt N\le  \frac98{\alpha\|x\|_1}.
 \end{eqnarray*}
It follows that the distance between $S$ and $T$  is at least
\[
\|x\|_1-2 \bigl(\tfrac98\alpha\|{x}\|_1 \bigr)\ge
\bigl(1-\tfrac94\alpha \bigr)\|{x}\|_1.
\]
So, by equation~(\ref{amasfini}), we have
\begin{eqnarray*}
&& \P_p\bigl(S\communique T,   |C(S)|<+\infty\bigr) \\
&&\qquad \le
\sum_{(3/8){C_d\alpha\|x\|_1}\le\|s\|_1\le (9/8){\alpha\|x\|_1}
}\P_p \bigl(s\communique \partial B \bigl(s, \bigl(1-\tfrac94\alpha \bigr)\|{x}\|_1 \bigr),   |C(s)|<+\infty \bigr)\\
&&\qquad \le   \bigl(1+2\times\tfrac98   \alpha\|{x}\|_1 \bigr)^{d}A_2e^{-B_2  (1-(9/4)\alpha )\|{x}\|_1}.
\end{eqnarray*}
So with a probability tending to $1$ exponentially fast with $\|{x}\|_1$, $S$ and $T$ belong to the infinite cluster.
Now,
\begin{eqnarray*}
&&\P_p \bigl(S\communique 0,   D(0,S)\ge {\eta}\|x\|_1\bigr) \\
&&\qquad \le  \sum_{(3/8){C_d\alpha\|x\|_1}\le\|s\|_1\le (9/8)\alpha\|x\|_1} \P_p \bigl(0\communique s,   D(0,s)\ge {\eta}\|x\|_1\bigr)\\
&&\qquad \le  \sum_{(3/8){C_d\alpha\|x\|_1}\le\|s\|_1\le (9/8)\alpha\|x\|_1} \P_p\bigl(0\communique s,   D(0,s)\ge \rho\|s\|_1\bigr)\\
&&\qquad \le  \sum_{(3/8){C_d\alpha\|x\|_1}\le\|s\|_1\le (9/8)\alpha\|x\|_1} A_1e^{-B_1\|s\|_1}\\
&&\qquad \le   \bigl(1+2\times\tfrac98   \alpha\|{x}\|_1 \bigr)^{d}A_1e^{-B_1(3/8)C_d{\alpha\|{x}\|_1}}.
\end{eqnarray*}
The second inequality is due to the choice (\ref{eta}) for $\eta$, the third to the result of Antal and Pisztora (\ref{antal_pisztora}).
We have
\begin{eqnarray*}
\|x-T\|_1 & \le & \bigl\|x-\|x\|_1\hat r\bigr\|_1+\bigl\|\|x\|_1\hat r-T\bigr\|_1 \\
& \le &  \biggl(\frac{C_d}3+\frac98 \biggr)\alpha\|x\|_1
 \le 2\alpha\|x\|_1= \frac{\eta}\rho \|x\|_1,
\end{eqnarray*}
and thus, similarly,
\begin{eqnarray*}
&& \P_p \bigl(T \communique x,   D(x,T)\ge {\eta}\|x\|_1\bigr) \\
&&\qquad \le  \P_p \bigl(T\communique x,   D(x,T)\ge \rho\|x-T\|_1\bigr)\\
&&\qquad \le  \sum_{(3/8){C_d\alpha\|x\|_1}\le\|t-\|x\|_1\hat{r}\|_1\le (9/8)\alpha\|x\|_1} \P_p\bigl(x\communique t,   D(x,t)\ge \rho\|x-t\|_1\bigr)\\
&&\qquad \le  \sum_{(3/8){C_d\alpha\|x\|_1}\le\|t-\|x\|_1\hat{r}\|_1\le (9/8)\alpha\|x\|_1} A_1e^{-B_1\|x-t\|_1}\\
&&\qquad \le   \bigl(1+2\times\tfrac98   \alpha\|{x}\|_1 \bigr)^{d}A_1e^{-B_1(1/{24})C_d{\alpha\|{x}\|_1}},
\end{eqnarray*}
where the last inequality follows from the estimate
\begin{eqnarray*}
\|x-t\|_1\ge \bigl\|\|x\|_1\hat r-t\bigr\|_1-\bigl\|x-\|x\|_1\hat r\bigr\|_1 & \ge &
\frac{C_d}{24}\alpha\|x\|_1.
\end{eqnarray*}

Denote by $\tilde{G}$ the event $G$ seen not as an event in the
macroscopic  percolation, but as a set of configurations of the microscopic percolation. Note that the event
\[
\tilde{G} \cap  \{ |C(S)|=+\infty  \}
\cap  \{ S\communique 0, D(0,S)\le {\eta}\|x\|_1  \}
\cap  \{ T\communique x, D(x,T)\le {\eta}\|x\|_1  \}
\]
is included in
\[
 \biggl\{ 0\communique x,
D(0,x) \leq   \biggl(
 \biggl( 1+\frac{3\eta}{2\rho}  \biggr)  ( 1+\eta  )^2  \mu(\hat x) +2\eta  \biggr)
\|{x}\|_1 \biggr\}.
\]
Thus, using equation (\ref{eta}) for $\eta$, and collecting all our previous estimates, we obtain
\begin{eqnarray*}
&&\P_p  \bigl(
0\communique x,
D(0,x)> \mu(\hat x)(1+ \epsilon) \|{x}\|_1  \bigr)\\
&&\qquad \le  \P_p  \biggl(
0\communique x,
D(0,x)>  \biggl(
 \biggl( 1+\frac{3\eta}{2\rho}  \biggr)  ( 1+\eta  )^2
\mu(\hat x)+2\eta  \biggr)
\|{x}\|_1  \biggr)\\
&&\qquad\le \P_p  ( \tilde{G}^c   )
+ \P_p  \bigl( S\communique T, |C(S)|<+\infty  \bigr)
+ \P_p  \bigl( S\communique 0, D(0,S)\ge {\eta}\|x\|_1  \bigr) \\
&&\qquad\quad{}+ \P_p  \bigl( T\communique x, D(x,T)\ge {\eta}\|x\|_1  \bigr)\\
&&\qquad\le 2Ae^{-({B}/{(NM)})\alpha\|{x}\|_1}+ \biggl(1+\frac{C\alpha}{NM}\|{x}\|_1
\biggr)^{2d}Ae^{-({B}/{(NM)})
 (1-2\alpha)\|{x}\|_1}\\
&&\qquad\quad{}
+  \biggl(1+2\times\frac98 \alpha\|{x}\|_1 \biggr)^{d}A_2e^{-B_2 (1-(9/4)\alpha)\|{x}\|_1}\\
&&\qquad\quad{}
+2  \biggl(1+2\times \frac98   \alpha\|{x}\|_1 \biggr)^{d}A_1e^{-B_1{\alpha\|{x}\|_1}/2},
\end{eqnarray*}
which ends the proof of the theorem.

\section{Lower large deviations}
\label{lower}
The aim of this section is to prove the lower large deviations estimate
for the chemical distance given by Theorem~\ref{devinf}.
First, we introduce some
definitions linked to the convexity of the asymptotic shape $\mathcal B_\mu
(0,1)$. Then, in Lemma~\ref{tempsplan}, we study the minimal number of open
edges needed to reach a given hyperplane at distance $r$ of the origin,
and, in Lemma~\ref{lesboites}, the minimal number of open
edges needed to cross a parallelepipedic
box. Finally, we prove the lower large deviations results.

\subsection{Definitions}
For each $y\in\Rd\backslash\{0\}$, the ball $\mathcal{B}_{\mu}(0,\mu(y))$
is a convex set, so
one can find a vector $n_y\in \mathcal{S}_2$ such that the linear form $\phi_y$
defined by
\[
\forall x\in\Rd\qquad \phi_y(x)= \langle n_y,x\rangle
\]
satisfies to $\phi_y(y)\ge 0$ and to
\begin{equation}
\label{separation}
\forall x \in \Rd \qquad  \bigl(  \mu(x)\le\mu(y) \bigr)  \quad
\Longrightarrow\quad
   \bigl(\phi_y(x)\le \phi_y(y)  \bigr).
\end{equation}

Note that the choice of $n_y$ is not necessarily unique. Using the fact that
the norms $\mu$ and $\|\cdot\|_2$ are homogeneous, it is possible to choose
the vector $n_y$ in such a
way that for each $y\in\Rd\backslash\{0\}$ and each $r>0$, one has
$n_{ry}=n_y$. In the following, we associate to every $y \in
\Rd\backslash\{0\}$ a unique $n_y$ satisfying these properties. We also
introduce the hyperplane $H_{y}=\ker \phi_y=(n_y)^{\perp}\!\dvtx$ geometrically
speaking, $y+H_{y}$ is a support hyperplane of the convex set $\mathcal{B}_{\mu}(0,\mu(y))$
at the point~$y$.

\begin{figure}

\includegraphics{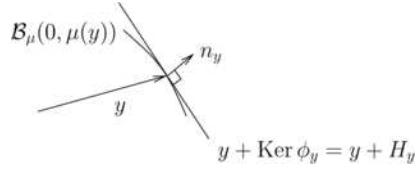}

\caption{Support hyperplane of the convex set
$\mathcal{B}_{\mu}(0,\mu(y))$ at the point~$y$\textup{:} the hyperplane $H_y$ may not be orthogonal to $y$.}
\label{f2}
\end{figure}

For $y\in\Rd\backslash\{0\}$ and $r \in \R^+ \backslash \{0\}$, note
\[
S_{y}^0=\{x\in\Rd\dvtx  \phi_{y}(x)<\phi_{y}(y)\}\quad\mbox{and}\quad S_{y}^\infty=\{x\in\Rd\dvtx  \phi_{y}(x)>\phi_{y}(y)\}.
\]
Then $S_{y}^0$---respectively, $S_{y}^\infty$---is the open half-space, delimited by the support hyperplane $y+H_{y}$ of  $\mathcal{B}_{\mu}(0,\mu(y))$ at the point~$y$,
containing---respectively, not containing---the origin (see Figure \ref{f2}).

The aim of the next lemma is to obtain, uniformly in $y \in \mathcal S_2$,
a bound on the norm of
points in the half-plane $S_{ry}^\infty$:

\begin{lemme}
\label{min_unif}
There exist two constants $c_d,c'_d>0$ such that
  for every $y \in \mathcal{S}_2$, we have:
\begin{itemize}
\item[$\bullet$] $\langle y,n_y\rangle \ge c_d$.
\item[$\bullet$] For every $r>0$,
$\displaystyle \inf \{ \|z\|_2,   z \in S_{ry}^\infty \} \ge c_d r$.
\item[$\bullet$] For every $r>0$,
$\displaystyle \inf \{ \|z\|_1,   z \in S_{ry}^\infty \} \ge c'_d r$.
\end{itemize}
\end{lemme}

\begin{pf}
As $\mu$ is a norm, it is equivalent to $\|\cdot\|_2$, so there exists
$K_1,K_2\in (0,+\infty)$ such that
\[
\forall x\in\Rd\qquad K_1\|x\|_2\le \mu(x)\le K_2\|x\|_2.
\]
Note $c_d=\frac{K_1}{K_2}$.
Choose $y \in \mathcal{S}_2$ and set $x=c_d n_y$.
We have $\mu(x)\le K_2\|x\|_2=K_2 c_d=K_1\le \mu(y)$.
It follows that $\phi_y(x)\le \phi_y(y)$, or equivalently
\[
c_d=\langle x,n_y\rangle\le \langle y,n_y\rangle.
\]
Now, let $z\in S_{ry}^\infty$: we have
\[
\|z\|_2\ge \langle z,n_{ry}\rangle=\phi_{ry}(z)>\phi_{ry}(ry)
=\langle ry,n_{ry}\rangle=r\langle y,n_y\rangle\ge c_d r.
\]
The last point is clear by the norm equivalence.
\end{pf}

\subsection{Passage-time from a point to a hyperplane}
Choose a direction $y \in \mathcal{S}_2$.
Define then for $r>0$
\[
b_{y}(r)=\inf\{D(0,z)\dvtx   z \in S_{ry}^\infty\}.
\]
This quantity is analogous to the usual passage time between the origin and a
hyperplane orthogonal to the first-coordinate axis at distance $r$ of the origin.
In this special case, $y=(1,0, \dots,0)$, and thanks to the symmetries of the grid, the direction of the support hyperplane of $\mathcal{B}_{\mu}(0,\mu(y))$ at the point $y$ is orthogonal to~$y$. But in a general direction $y$, the relevant hyperplane for the growth of the set~$B_t$ of wet vertices at time $t$ in the direction $y$ is
$H_{y }$, which does not need to be orthogonal to $y$ (see Figure \ref{f2}).
As in the paper by Cox and Durrett~\cite{cd81}, we can study this quantity by
using the
asymptotic shape result given in \cite{gm04}:
for every $p>p_c(d)$,
\begin{equation}
\label{asymptotic-shape}
 \lim_{t\to +\infty}\mathcal{D} \biggl(\frac{B_t}t,\mathcal{B}_{\mu}(0,1) \biggr)=
0,\qquad\Pcond_p\as
\end{equation}
and obtain the following lemma:

\begin{lemme}
\label{tempsplan}
\[
\displaystyle \sup_{y \in \mathcal{S}_2}  \biggl| \frac{b_{y}(r)}r -\mu(y)  \biggr| \to 0,  \qquad   \Pcond_p\as
\]
\end{lemme}

\begin{pf}
As we
work under $\Pcond_p$, we restrict ourselves to the event $\{0 \communique
\infty\}$.

Let $\epsilon>0$. By the asymptotic
shape result (\ref{asymptotic-shape}),
there $\Pcond_p\as$ exists a random $T$ such that
\begin{equation}
\label{blabla}
\forall t \ge T\qquad\frac{B_t}t \subset (1+\epsilon)
\mathcal{B}_{\mu}\quad\mbox{and}\quad\mathcal{B}_{\mu} \subset
\frac{B_t}t+\mathcal{B}_{\mu}(0, \epsilon).
\end{equation}
For every $r>0$, for every $y \in \mathcal{S}_2$,
there exists a point $z_r^{y} \in  S_{ry}^\infty$ such that
$b_{y}(r)=D(0,z_r^{y}) \ge \|z_r^y\|_1$.
Note that, by the previous lemma, we have  $\|z_r^y\|_1 \ge \inf
\{\|z\|_1\dvtx
  z
\in S_{ry}^\infty\}\ge c_d' r$,
and thus, as soon as $c'_d r \ge T$, we have
\[
\forall y \in \mathcal{S}_2 \qquad \frac{z_r^{y}}{ b_{y}(r)} \in \frac{B_{b_{y}(r)}}{b_{y}(r)}
\subset (1+\epsilon)\mathcal{B}_{\mu}.
\]
This implies, by definition of $H_{y}$ and by convexity of
$\mathcal{B}_{\mu}$,  that $\frac1{b_{y}(r)}\frac{\langle z_r^{y}, n_{y} \rangle  }{\langle y, n_{y} \rangle }y \in
(1+\epsilon)\mathcal{B}_{\mu}$ and thus
$\Pcond_p\as$, for all $r$ large enough,
\[
\label{controle-b}
\forall y \in \mathcal{S}_2 \qquad \frac{r}{b_{y}(r)} \le \frac{1+\epsilon}{\mu({y})}.
\]

On the other hand, by definition of $b_{y}(r)$, we have $S_{ry}^\infty \cap B_{b_{y}(r)-1} = \varnothing$, or, in other words,
$B_{b_{y}(r)-1} \subset S_{ry}^0$. By dilatation, we obtain
\[
\frac{B_{b_{y}(r)-1}}{b_{y}(r)-1}
\subset S_{{ry}/({b_{y}(r)-1})}^0,
\]
and by definition of $H_{y}$, this leads to:
\[
\frac{B_{b_{y}(r)-1}}{b_{y}(r)-1}+\mathcal{B}_{\mu}(0, \epsilon)
\subset S_{( {r}/({b_{y}(r)-1})+\epsilon)y}^0.
\]
Using (\ref{blabla}),  we obtain, as soon as $c'_d r -1\ge T$,
\[
\forall y \in \mathcal{S}_2 \qquad \mathcal{B}_{\mu} \subset \frac{B_{b_{y}(r)-1}}{b_{y}(r)-1}
+\mathcal{B}_{\mu}(0, \epsilon),
\]
 and thus $\mathcal{B}_{\mu} \cap  S_{({r}/({b_{y}(r)-1})+\epsilon)y}^\infty=\varnothing$, leading to
\[
\mu\biggl ( \biggl(\frac{r}{b_{y}(r)-1}+\epsilon \biggr)y \biggr)
=    \biggl( \frac{r}{b_{y}(r)-1} + \epsilon  \biggr) \mu(y)
>1.
\]
Finally, we get $\Pcond_p\as$, for all $r$ large
enough,
\[
\forall y \in \mathcal{S}_2 \qquad \frac{r}{b_{y}(r)-1} \ge
\frac{1}{\mu({y})}-\epsilon,
\]
which ends the proof.
\end{pf}

\subsection{Crossings of parallelepipedic boxes}
We want first to find $d$ directions $(y_1=y, y_2, \dots ,
  y_d)$ with $y_i \in \mathcal{S}_2$ such that the asymptotic time constants are the same along all these directions and such that the directions of the support hyperplanes of $\mathcal{B}_{\mu}$ in these directions are linearly independent.

For $y \in \mathcal{S}_2$, consider the vector $n_{y} \in \mathcal{S}_2$ orthogonal to a support hyperplane
of~$\mathcal{B}_{\mu}$ in the direction $y$ as defined previously, and the isometries $(g_{n_y,2},
\dots,$ $  g_{n_y,d}) \in  ( \mathcal O  ( \Zd  )  )^{d-1}$  given by equation~(\ref{controle-normales}) in Lemma~\ref{controle-direction} and set
\begin{eqnarray*}
(n_1, n_2, \dots ,
  n_d) & = & \bigl( n_y, g_{n_y,2}( n_y), \dots,  g_{n_y,d}( n_y)\bigr), \\
(y_1, y_2, \dots ,
  y_d) & = & \bigl(y, g_{n_y,2}(y), \dots,  g_{n_y,d}(y)\bigr).
\end{eqnarray*}

For $k \in \Zd$, $\alpha \in (\R_+^*)^d$, we define boxes adapted to study the travel times in the directions $y_1, y_2, \dots ,
  y_d$; they are analogous to the rectangular boxes introduced to estimate the travel time in the first-coordinate axis in classical first-passage percolation.
\begin{eqnarray*}
T_{(y)}(k,\alpha) & = &  \biggl\{v \in \Zd\dvtx \, \forall m \in \{1,
\dots,d\}\,
k_m \le \frac{\langle v,n_m \rangle}{\langle y_m,n_m \rangle} <k_m+\alpha_m \biggr\},\\
\partial_-^m T_{(y)}(k,\alpha) & = &
 \left\{
 \begin{array}{l}
v \in \Zd \backslash T_{(y)}(k,\alpha)\dvtx  \cr
      \quad\bullet \forall j \in \{1, \dots,d\} \backslash \{m\} \,
 k_j \le \dfrac{\langle v,n_j \rangle}{\langle y_j,n_j \rangle} <k_j+\alpha_j
 \cr
     \quad \bullet \dfrac{\langle v,n_m \rangle}{\langle y_m,n_m \rangle}
     <k_m\cr
     \quad \bullet \exists w \in T_{(y)}(k,\alpha) \, \|w - v\|_1=1
\end{array}
 \right\},
\\
\partial_+^m T_{(y)}(k,\alpha) & = &
 \left\{
\begin{array}{l}
v \in \Zd \backslash T_{(y)}(k,\alpha):  \cr
      \quad\bullet \forall j \in \{1, \dots,d\} \backslash \{m\} \quad
 k_j \le \dfrac{\langle v,n_j \rangle}{\langle y_j,n_j \rangle} <k_j+\alpha_j
 \cr
      \quad\bullet \dfrac{\langle v,n_m \rangle}{\langle y_m,n_m \rangle} \ge
      k_m+\alpha_m\cr
     \quad \bullet \exists w \in T_{(y)}(k,\alpha) \, \|w - v\|_1=1
\end{array}
 \right\}.
\end{eqnarray*}
We can now define, using the same terminology as in first-passage
percolation,  the crossing time of the box $T_{(y)}(k,\alpha)$ in the
$m$th direction:
\[
t^m_{(y)}(k,\alpha)  =  \inf  \left\{
\begin{array}{l}
|\gamma|,\mbox{ where $\gamma$ is an open path from} \cr
\mbox{a point in $\partial_-^m T_{(y)}(k,\alpha)$
to a point in $\partial_+^m T_{(y)}(k,\alpha)$}\cr
\mbox{included but its ends in $T_{(y)}(k,\alpha)$}
\end{array}
 \right\}.
\]
The next lemma gives a convergence in probability, uniformly in the
direction $y$, of these minimal crossing times of boxes:

\begin{lemme}
\label{lesboites}
Let $\alpha=(\alpha_1, \dots, \alpha_d) \in (0,+\infty)^d$.
Then for every $\epsilon>0$ we have
\[
\forall m\in \{1,\dots,d\} \qquad
\lim_{r\to +\infty}\sup_{k \in \Zd} \sup_{y \in\mathcal{S}_2} \P_p \bigl(t^m_{(y)}(k,r\alpha) \leq \bigl(\mu(y)-\epsilon\bigr)r\alpha_m\bigr)=0.
\]
\end{lemme}

\begin{pf}
Fix $\alpha=(\alpha_1, \dots, \alpha_d) \in (0,+\infty)^d$ and $\epsilon>0$. Choose
\begin{equation}
\label{eta2}
\eta= \frac{C_d \epsilon}{6\rho d^{5/2}}  \biggl(\min_{1 \le m \le d }\alpha_m \biggr)>0,
\end{equation}
where $C_d$ is given in Lemma~\ref{controle-direction} and $\rho$ is the
constant introduced by Antal and Pisztora---see equation (\ref{antal_pisztora}); set $\tau=\eta r$.

Choose $k \in \Zd$ and $y\in\mathcal{S}_2 $.

We introduce the
following partition of $\Zd$ into boxes of size $\tau$, which tends to
infinity when $r$ goes to infinity, but will still be small when compared to $r$: for every $\overline x\in \Zd$, and every $A \subset \Zd$
\begin{eqnarray*}
\boite_{(y)}(\overline x)& = &  \biggl\{v \in \Zd\dvtx \, \forall m \in \{1,
  \dots,d\}\,
\overline x_m \tau \le \frac{\langle v,n_m \rangle}{\langle y_m,n_m
  \rangle} <(\overline x_m+1) \tau \biggr\},\\
\voisinage_{(y)} ( A  )& = & \bigl\{\overline x \in \Zd\dvtx \,
 \boite_{(y)}(\overline x) \cap A \neq \varnothing \bigr \}.
\end{eqnarray*}
Then $\voisinage_{(y)}  ( \partial_-^m T_{(y)}(k,r\alpha)  )$, which plays the role of an approximation at a larger scale of $\partial_-^m T_{(y)}(k,r\alpha)$, is
$*$-connected and a simple estimation leads to
\[
\bigl| \voisinage_{(y)}  \bigl( \partial_-^m T_{(y)}(k,r\alpha)  \bigr)\bigr|
\le  \biggl(2 +\frac{1}{\eta r}  \biggr) \mathop{\prod_{1\le j\le d}}_{j \neq m}  \biggl(2+
  \frac{\alpha_j}{\eta}  \biggr).
\]
The partition into boxes with size $\tau$ was introduced to obtain this
estimate: while $| \partial_-^m T_{(y)}(k,r\alpha)|$ is of order $r^{d-1}$, the cardinal $| \voisinage_{(y)}  ( \partial_-^m T_{(y)}(k,r\alpha)  )|$ of its approximation with large boxes remains bounded when $r$ goes to infinity.
If $v \in \partial_-^m T_{(y)}(k,r\alpha)$, then there exists a unique $\overline x_v
\in \voisinage_{(y)}  ( \partial_-^m T_{(y)}(k,r\alpha)  )$ such that $v
\in \boite_{(y)}(\overline x_v)$. We define
\[
\mathcal W(v)=\bigl\{\overline x \in \voisinage_{(y)}  \bigl(
  \partial_-^m T_{(y)}(k,r\alpha)  \bigr)\dvtx   \|\overline x-\overline x_v\|_\infty=2\bigr\}\neq
\varnothing.
\]
If $\boite_{(y)}(\overline x) \cap C_\infty \neq \varnothing$, then define $c(\overline x)$
as the point in $\boite_{(y)}(\overline x) \cap C_\infty$ which is the closest to
$\tau(x+(1/2,\dots,1/2))$ (use the lexicographic order if necessary). Remember that the box $\boite_{(y)}(\overline x)$ has size $\tau$, which tends to infinity with $r$, and thus as $r$ goes to infinity, we expect the probability that $\boite_{(y)}(\overline x) \cap C_\infty \neq \varnothing$ to go to $1$.

Now, in the following inequality, we approximate the event $\{t^m_{(y)}(k,r\alpha) \leq (\mu(y)-\epsilon)r\alpha_m\} $ by the event in (\ref{AA}), and the three last terms correspond to the difference between them, and are expected to be
small:
\begin{eqnarray}
&&\P_p \bigl(t^m_{(y)}(k,r\alpha) \leq \bigl(\mu(y)-\epsilon\bigr)r\alpha_m\bigr)
\nonumber\label{AA}\\[-8pt]\\[-8pt]
&&\qquad \leq  \P_p  \left(\begin{array}{l}
\bullet \forall \overline x \in \voisinage_{(y)}  \bigl( \partial_-^m T_{(y)}(k,r\alpha)
 \bigr) \, \boite_{(y)}(\overline x) \cap C_\infty \neq \varnothing \nonumber\\
\bullet \forall v \in \partial_-^m T_{(y)}(k,r\alpha) \, \bigl(v \communique
\partial_+^m T_{(y)}(k,r\alpha) \Rightarrow v \communique \infty\bigr) \\
\bullet \exists v \in \partial_-^m T_{(y)}(k,r\alpha) \mbox{ such that } \\
      \quad \exists w \in
\partial_+^m T_{(y)}(k,r\alpha) \, D(v,w) \leq \bigl(\mu(y)-\epsilon\bigr)r\alpha_m
 \\
      \quad \exists \overline x \in \mathcal W(v) \, D\bigl(v, c(\overline x)\bigr) < \rho \|v-c(\overline x)\|_1
\end{array}
 \right)\\
&&\qquad\quad{} + \P_p  \bigl( \exists \overline x \in \voisinage_{(y)} \bigl ( \partial_-^m T_{(y)}(k,r\alpha)
 \bigr) \, \boite_{(y)}(\overline x) \cap C_\infty = \varnothing  \bigr) \label{BB}\\
&&\qquad\quad{} + \P_p  \bigl( \exists v \in \partial_-^m T_{(y)}(k,r\alpha) \, v
  \communique \partial_+^m T_{(y)}(k,r\alpha), \, v \not\communique \infty
 \bigr)\label{CC}\\
&&\qquad\quad{} + \P_p \left (
\begin{array}{l}
\bullet \forall \overline x \in \voisinage_{(y)}  \bigl( \partial_-^m T_{(y)}(k,r\alpha)
 \bigr) \, \boite_{(y)}(\overline x) \cap C_\infty \neq \varnothing \\
\bullet \exists v \in \partial_-^m T_{(y)}(k,r\alpha) \, \forall  \overline x \in \mathcal W(v) \\
       \quad D\bigl(v, c(\overline x)\bigr) \ge \rho \|v-c(\overline x)\|_1
\end{array}
 \right)\label{DD}.
\end{eqnarray}
Let us estimate the three error terms first.

\textit{Estimate for }(\ref{BB}).
Let $\overline x \in \Zd$, and let us prove that $\boite_{(y)}(\overline x)$
contains a ball for the norm $\|\cdot\|_2$ with radius proportional to $\tau$. Let us introduce first the point $a_{\overline x} \in \Rd$, which represents the center of $\boite_{(y)}(\overline x)$, such that
\[
\forall m \in \{1, \dots, d\} \qquad
 \frac{\langle a_{\overline x},n_m\rangle}{\langle y_m,n_m\rangle}= \biggl( \overline x_m +\frac12  \biggr) \tau.
\]
Then we have, with $c_d$ given by Lemma~\ref{min_unif},
\begin{eqnarray*}
\boite_{(y)}(\overline x) & = &  \biggl\{v \in \Zd\dvtx \quad \forall m \in \{1, \dots, d\}
\,
-\frac{\tau}2 \le \frac{\langle \overline x-a_{\overline x},n_m\rangle}{\langle y_m,n_m\rangle}<
\frac{\tau}2  \biggr\} \\
& \supset &  \biggl\{v \in \Zd\dvtx \quad \|\overline x-a_{\overline x}\|_2 \le \frac{\tau c_d}2 \biggr\}.
\end{eqnarray*}
Thus the box $\boite_{(y)}(\overline x)$ contains the ball $\mathcal{B}_2(a_{\overline x},\tau c_d/2)$; this radius does not depend on the direction $y$.

Using then equation (\ref{amasinfini}), we get
\begin{eqnarray*}
\mbox{(\ref{BB})} & \leq  & \bigl|\voisinage_{(y)}  \bigl( \partial_-^m T_{(y)}(k,r\alpha)
   \bigr)\bigr| \sup_{z \in \Zd}\P_p  \bigl( \mathcal{B}_2(z,\tau c_d/2)\cap C_\infty = \varnothing  \bigr)
   \\
 & \leq  &  \biggl(2 +\frac{1}{\eta r}  \biggr) \mathop{\prod_{1\le j\le d}}_{j \neq m}  \biggl(2+
  \frac{\alpha_j}{\eta}  \biggr) A_3 \exp  \biggl(-\frac{B_3 c_d \eta r}{2\sqrt d} \biggr) ,
\end{eqnarray*}
which tends to $0$ when $r$ goes to infinity, uniformly in the direction $y \in \mathcal{S}_{2}$.

\textit{Estimate for }(\ref{CC}).
Note that if $v \in \partial_-^m T_{(y)}(k,r\alpha)$ and $w  \in \partial_+^m T_{(y)}(k,r\alpha)$, then, by definition of the box, we have
\[
\biggl|\frac{\langle w-v,n_m \rangle}{\langle y_m,n_m\rangle}  \biggr|\ge r \alpha_m.
\]
By Lemma~\ref{min_unif}, we have
\[
\sqrt d \|w-v\|_1 \ge \|w-v\|_2 \ge| \langle w-v,n_m \rangle |\ge c_d r \alpha_m,
\]
and thus, using translation invariance and equation (\ref{amasfini}),
\begin{eqnarray*}
\mbox{(\ref{CC})} & \leq  & \bigl|\partial_-^m T_{(y)}(k,r\alpha)\bigr| \P_p  \biggl( |C(0)| \geq
\frac{c_d}{\sqrt d} r  \alpha_m ,   0 \not\communique \infty  \biggr) \\
 & \leq  & K_\alpha r^{d-1} A_2 \exp  \biggl(-B_2
\frac{c_d}{\sqrt d} r  \alpha_m  \biggr),
\end{eqnarray*}
which tends to $0$ when $r$ goes to infinity, uniformly in the direction
$y \in \mathcal{S}_{2}$. ($K_\alpha$ is a constant depending only on the
dimension and $\alpha$.)

\textit{Estimate for }(\ref{DD}).
Let $v \in \partial_-^m T_{(y)}(k,r\alpha)$ and $\overline x \in \mathcal W(v)$
such that $  \boite_{(y)}(\overline x) \cap C_\infty \neq \varnothing$.
Then by construction,
\[
\max_{1 \le m \le d}  \biggl|\frac{\langle v-c(\overline x),n_m \rangle}{\langle y_m,n_m\rangle}\biggr |
 \ge \tau.
 \]
But then, using Lemma~\ref{min_unif}, we obtain
\begin{eqnarray*}
\|v-c(\overline x)\|_1 & \ge &\frac1{\sqrt d} \|v-c(\overline x)\|_2
\ge  \frac1{\sqrt d} \max_{1 \le m \le d}  | \langle v-c(\overline x),n_m \rangle  | \\
& \ge  & \frac{c_d}{\sqrt d} \max_{1 \le m \le d}  \biggl|\frac{\langle v-c(\overline x),n_m \rangle}{\langle y_m,n_m\rangle} \biggr|
\ge  \frac{c_d\tau}{\sqrt d}.
\end{eqnarray*}
Using Antal and Pisztora's result (\ref{antal_pisztora}), one has the following bound:
\begin{eqnarray*}
(\ref{DD}) & \leq  & \bigl|\partial_-^m T_{(y)}(k,r\alpha)\bigr| \times 5^{d-1}2d
A_1 \exp \biggl ( -B_1 \frac{c_d\eta r}{\sqrt d}  \biggr),
\end{eqnarray*}
which tends to $0$ when $r$ goes to infinity, uniformly in the direction
$y \in \mathcal{S}_{2}$.

\textit{Estimate for }(\ref{AA}).
Let $v \in \partial_-^m T_{(y)}(k,r\alpha)$ and $\overline x \in \mathcal W(v)$ such that $  \boite_{(y)}(\overline x) \cap C_\infty \neq \varnothing$. Then by construction,
\[
\max_{1 \le m \le d}  \biggl|\frac{\langle v-c(\overline x),n_m
\rangle}{\langle y_m,n_m\rangle} \biggr| \le 3\tau.
\]
But then, using equation~(\ref{controle-normales}) in Lemma~\ref{controle-direction}, we have
\begin{eqnarray*}
\|v-c(\overline x)\|_1 & \le & \sqrt d \|v-c(\overline x)\|_2
\le \frac{d^2}{C_d} \|v-c(\overline x)\|_{(n_i),2} \\
& \le & \frac{d^{5/2}}{C_d} \max_{1 \le m \le d}  |\langle v-c(\overline x),n_m \rangle
|\\
&\le& \frac{d^{5/2}}{C_d} \max_{1 \le m \le d}  \biggl|\frac{\langle v-c(\overline x),n_m \rangle}{\langle y_m,n_m\rangle} \biggr| \\
& \le & \frac{3d^{5/2}\tau}{C_d}.
\end{eqnarray*}
We obtain
\begin{eqnarray*}
(\ref{AA}) & \leq  &
 \P_p \left (
\begin{array}{l}
\exists \overline x \in \voisinage_{(y)}  \bigl( \partial_-^m T_{(y)}(k,r\alpha)  \bigr)
\mbox{ such
    that } \\
  \quad \bullet \boite_{(y)}(\overline x ) \cap C_\infty \neq \varnothing, \\
  \quad \bullet \exists w \in \partial_+^m T_{(y)}(k,r\alpha) \\
     \qquad D\bigl(c(\overline x),w\bigr) < \dfrac{3\rho d^{5/2}\eta r}{C_d} + \bigl(\mu( y)-\epsilon\bigr)r\alpha_m
\end{array}
 \right).
\end{eqnarray*}
By the choice (\ref{eta2}) we made for $\eta$, we obtain that
\begin{eqnarray*}
(\ref{AA}) & \leq  & \bigl|\voisinage_{(y)}  \bigl( \partial_-^m T_{(y)}(k,r\alpha)
\bigr)\bigr|
\\
&&{} \times \sup_{\overline x \in \Zd} \P_p  \left(
\begin{array}{l}
\bullet \boite_{(y)}(\overline x) \cap C_\infty \neq \varnothing, \\
\bullet \exists y \in \partial_+^m T_{(y)}(k,r\alpha),\,
D\bigl(c(\overline x),y\bigr) \leq \biggl(\mu(y)-\dfrac{\epsilon}2\biggr)r\alpha_m
\end{array}
\right ) \\
& \leq &  \biggl(2 +\frac{1}{\eta r}  \biggr) \mathop{\prod_{1\le j\le d}}_{j \neq m} \biggl (2+
  \frac{\alpha_j}{\eta}  \biggr) \P_p  \biggl( b_{y_m}(\alpha_m r) \leq
 \biggl( \mu(y)-\frac{\epsilon}2  \biggr)\alpha_m r  \biggr)
\end{eqnarray*}
which tends to $0$ when $r$ goes to $\infty$ by Lemma~\ref{tempsplan}. Note that this convergence is uniform in $k \in \Zd$ and $y \in \mathcal{S}_2$.
\end{pf}

\subsection[Lower large deviations\textup{:} proof of Theorem~1.6]{Lower large deviations\textup{:} proof of Theorem~\textup{\protect\ref{devinf}}}
We essentially follow the main lines of the proof in the classical case
 by Grimmett and Kesten~\cite{gk84}:
A ``too short'' path should cross ``many'' boxes in a ``too short'' time, and by the previous result and a counting argument, this probability can be made exponentially small. The two main difficulties are to deal with geometric problems due to the
fact that we want large deviations not only along the coordinate axes, but
in all directions and the uniformity we require in this direction.



\textit{Step  }1. Definition of boxes adapted to direction $y$.

Choose $M$ and $N$ large enough, that will be fixed later.

For $\overline k=(\overline k_1,\dots,\overline k_d) \in \Zd$, we define
\begin{eqnarray*}
S_{(y)}(\overline k) & = &  \biggl\{v \in \Zd\dvtx \, \forall m \in \{1, \dots,d\}  \, N \overline k_m \le
\frac{\langle v,n_m \rangle}{\langle y_m,n_m\rangle}<N(\overline k_m+1)\biggr \},\\
T_{(y)}(\overline k) & = &  \left\{
\begin{array}{l}
v \in \Zd\dvtx   \forall m \in \{1, \dots,d\}   \\
\quad N \overline k_m -M\le
\dfrac{\langle v,n_m \rangle}{\langle y_m,n_m\rangle} <N(\overline k_m+1) +M
\end{array}
 \right\}.
\end{eqnarray*}
The $S_{(y)}(\overline k)$'s are large ``twisted square''  boxes, adapted to the studied direction of progression $y$ and its conjugates $(y_2, \dots,y_d)$, that induce a partition of
$\Zd$, and the $T_{(y)}(\overline k)$'s are still much larger boxes centered in the
$S_{(y)}(\overline k)$'s.

In $T_{(y)}(\overline k)$, the small box $S_{(y)}(\overline k)$ is
surrounded by $2d$ boxes of the type ($1 \leq m \leq d)$:
\begin{eqnarray*}
B_{(y),m}^+(\overline k)& =  &
 \left\{
\begin{array}{l}
v \in \Zd\dvtx \\
\quad  \forall j \neq m \,
N \overline k_j-M \le \dfrac{\langle v,n_j \rangle}{\langle y_j,n_j\rangle} <N(\overline k_j+1)+M, \\
\quad N (\overline k_m+1) \le \dfrac{\langle v,n_m \rangle}{\langle y_m,n_m\rangle} <N(\overline k_m+1)+M
\end{array}
 \right\} \\
B_{(y),m}^-(\overline k) & = &
 \left\{
\begin{array}{l}
v \in \Zd\dvtx \\
\quad  \forall j \neq m \,
N \overline k_j-M \le \dfrac{\langle v,n_j \rangle}{\langle y_j,n_j\rangle} <N(\overline k_j+1)+M, \\
\quad N \overline k_m-M \le \dfrac{\langle v,n_m \rangle}{\langle y_m,n_m\rangle} <N\overline k_m
\end{array}
\right \}.
\end{eqnarray*}
We define the inside border and the outside border of the box $B_{(y),m}^-(\overline k)$,
relatively to $T_{(y)}(\overline k)$:
\begin{eqnarray*}
\partial_{\mathrm{in}}B_{(y),m}^-(\overline k) & = & \bigl\{  v
\in T(\overline k)\backslash B_{(y),m}^-(\overline k)\dvtx \exists w
\in B_{(y),m}^-(\overline k) \, \|v -w\|_1=1 \bigr\},\\
\partial_{\mathrm{out}}B_{(y),m}^-(\overline k) & = &  \left\{
\begin{array}{l}
v \in T_{(y)}(\overline k_1,\dots, \overline k_{m-1},\overline k_m-1,\overline k_{m+1},
\dots, \overline k_d)\backslash T_{(y)}(\overline k)\dvtx \\
\quad\exists w
\in B_{(y),m}^-(\overline k) \, \|v -w\|_1=1
\end{array}  \right\},
\end{eqnarray*}
and the borders of the other boxes can be defined in the obvious analogous
manner.

The point is that a path, visiting $S_{(y)}(\overline k)$ and exiting from
$T_{(y)}(\overline k)$, has to cross one of these $2d$ boxes surrounding $S_{(y)}(\overline k)$ in $T_{(y)}(\overline k)$
from a point in its inside border to a point in its outside border, say $B_{(y),m}^+(\overline k)$, for instance. And,
roughly speaking, the fastest way to cross it is to follow the $y_m$
direction,  and this should take an amount of steps of order
$\mu(y_m )M=\mu(y)M$ if $M$ is large. We can easily estimate the size of borders of boxes:

\begin{lemme}
\label{bords}
There exist a strictly positive constant $K_d$  depending only
on the dimension $d$ and not on $y \in \mathcal{S}_2 $ such that for every $k
\in \Zd$,
for every $m \in \{1, \dots,d\}$,
\begin{eqnarray*}
\frac1{K_d} (N+2M)^{d-1} \le & \bigl|\partial_{\mathrm{in}}B_{(y),m}^+(k)\bigr| & \le K_d (N+2M)^{d-1} ,\\
\frac1{K_d}  (N+2M)^{d-1}  \le & \bigl|\partial_{\mathrm{out}}B_{(y),m}^+(k)\bigr| & \le K_d (N+2M)^{d-1}.
\end{eqnarray*}
The same is also true for the borders of $B_{(y),m}^-(k)$\textup{'}s.
\end{lemme}

\textit{Step  }2. Construction of crossings.

The construction is exactly the same as the one in
\cite{gk84}. We thus only give the way to adapt it.

Let $r>0$ large enough and let $\gamma=(v(0),\dots, v(\nu))$ be a path
from $0$ to a point in $S_{r y}^\infty$. We associate to $\gamma$
the following two sequences. First set $k(0)=0$ and $a(0)=0$. Let then $v(a(1))$ be the
first vertex along $\gamma$ to be outside $T_{(y)}(k(0))$, and let $k(1)$
be the
coordinates of the small box of type $S$ containing $v(a(1))$: $v(a(1)) \in S_{(y)}(k(1))$,
and build the two sequences recursively, to obtain $(a(1), \dots, a(\tau(\gamma)))$
and $(k(0),\dots, k(\tau(\gamma)))$ such that:
\begin{longlist}[3.]
\item[1.]
$0=a(0)<a(1)< \cdots <a(\tau(\gamma))\le \nu$,
\item[2.]
$v(a(i)) \in S_{(y)}(k(i))$,
\item[3.]
$a(i+1)$ is the smallest integer $a$ larger than $a(i)$ such that $v(a)
\notin T_{(y)}(k(i))$.
\end{longlist}
The final terms satisfy
\[
\forall a(\tau(\gamma)) \le j \le \nu \qquad v(j) \in T_{(y)}(k(\tau(\gamma))).
\]
Note that the portion $\gamma(i)$ of $\gamma$ between $v(a(i-1))$ and $v(a(i))$
has to cross one of the $2d$ boxes of type $B_{(y),m}(k(i-1))$
surrounding $S_{(y)}(k(i-1))$ in $T_{(y)}(k(i-1))$
from a point in its inside border to a point in its outside border.
 We are interested in these crossings and the amount of steps they use.

But first, by the process of loop removal described in \cite{gk84}, remove
the double points from $\Gamma=(k(0),\dots, k(\tau(\gamma)))$, and obtain
\[
\tilde \Gamma=\bigl(l(0),\dots, l(\sigma(\gamma))\bigr),
\]
where $l(a)=k(j_a)$ and $0<j_0<\cdots<j_{\sigma(\gamma)}\leq \tau(\gamma)$.
Note that although
we can have $j_{\sigma(\gamma)}< \tau(\gamma)$,
it is always true that $k(j_{\sigma(\gamma)})=k(\tau(\gamma))$. By construction,
\[
\forall m \in \{1, \dots, d\} \mbox{ } \forall j \in \{0, \dots,
\tau(\gamma)-1\} \qquad
 \biggl| \frac{\langle n_{y_m},k(j+1)-k(j) \rangle}{\langle n_{y_m},y_m \rangle}\biggr | \leq
\frac{M}N +1,
\]
and this property is preserved by the loop removal process in the following
sense:
\[
\forall m \in \{1, \dots, d\} \mbox{ } \forall j \in \{0, \dots, \sigma(\gamma)-1\} \qquad
 \biggl| \frac{\langle n_{y_m},l(j+1)-l(j) \rangle}{\langle n_{y_m},y_m \rangle} \biggr|  \leq
\frac{M}N +1.
\]

\textit{Step  }3. Coloring of crossings.

Consider the portion $\gamma(i)$ of $\gamma$ between $v(a(i-1))$ and
$v(a(i))$, and
define:
\[
\mathcal L(i)=\max_{m \in \{1, \dots,d\}}
 \biggl| \frac{\langle  n_{y_m},v(a(i)) -v(a(i-1)) \rangle}{\langle n_{y_m},y_m \rangle} \biggr |.
 \]
By construction, $M \le \mathcal L(i) \le M+N$ for $1 \le i \le \tau(i)$. Now, for $i
\in \{1, \dots, \sigma(\gamma)\}$, consider the portion $\gamma(j_i)$ between the two boxes
$S_{(y)}(k(j_i-1))$ and $S_{(y)}(k(j_i))$ and give to the vector $l(i)=k(j_i)$ the
color white if
\[
\gamma(j_i) \mbox{ is open}\quad\mbox{and}\quad |\gamma(j_i)| \leq \bigl(\mu(y)
-2\epsilon\bigr) \mathcal L(j_i),
\]
and in black otherwise. Denote by $w(\gamma)$ the number of white points in the
sequence $(l(1),\dots, l(\sigma(\gamma)))$ of crossings associated to $\gamma$.
The next lemma corresponds to Lemma 3.5
in the paper by Grimmett and Kesten~\cite{gk84}:

\begin{lemme}
Suppose that $\epsilon$, $r$, $M$ and  $N$  satisfy the following:
\begin{longlist}[3.]
\item[1.]
$\epsilon$ is small\textup{:} $ 0 <\epsilon <\min\{\mu( z)\dvtx   z \in \mathcal{S}_2\}$,
\item[2.]
$M/N$ is large\textup{:} $ \forall z \in \mathcal{S}_2 \mbox{ } M(\mu(z)-3\epsilon) \geq
(M+N) (\mu(z)-4\epsilon)$,
\item[3.]
$r$ is large\textup{:} $\forall z \in \mathcal{S}_2 \mbox{ } r \epsilon \ge
(M+2N)(\mu(z)-4\epsilon)$.
\end{longlist}
Then for any $y \in \mathcal S_2$, an open path $\gamma$, traveling from $0$ to $S_{r y}^\infty$  and whose length~$|\gamma|$ is less or equal to  $r (\mu(
y)-5\epsilon)$, satisfies
\[
w(\gamma) \ge \frac{\epsilon \sigma}{2\mu(e_1)}\quad \mbox{and}\quad
\sigma (\gamma)\ge \frac{r}{M+N}-1.
\]
\end{lemme}

\begin{pf}
The proof for a fixed $y \in \mathcal{S}_2$ is  exactly as in \cite{gk84}.
The only difference is the uniformity in $y \in \mathcal{S}_2$, that can be obtained  because $\mu$ is a norm and is bounded
away from $0$ and bounded away from infinity on the compact set $\mathcal{S}_2$.
\end{pf}

Note that the event $\{l(i) \mbox{ is white}\}$ in contained in the event
\[
\mathcal E_{(y),i} = \left\{
\begin{array}{l}
\mbox{a vertex in $S_{(y)}(k(j_i))$ is joined to a vertex } \\
\mbox{outside } T_{(y)}(k(j_i)) \mbox{ by a path using} \\
\mbox{less than }(M+N) (\mu(y)-2 \epsilon) \mbox{ steps}
\end{array}
 \right\}.
 \]
To conclude the proof by a counting argument, as these event are only locally dependent, it only remains to prove the following lemma:

\begin{lemme}
Let $0 <\epsilon<\mu_{\mathrm{inf}}=\min\{\mu( z)\dvtx   z \in \mathcal{S}_1\}$. Then
\[
\begin{array}{l}
p(M,N,\epsilon)  =
\mathop{\sup\limits_{k \in \Zd}}\limits_{y\in\mathcal{S}_2}  \P_p \left (
\begin{array}{l}
\exists v \in S_{(y)}(k) ,\, \exists w \notin T_{(y)}(k), \\
\quad v \communique w,   D(v,w) \leq  (M+N) \bigl(\mu(y)-2 \epsilon\bigr)
\end{array}
\right )
\end{array}
\]
goes to $0$ when $M$ and $N$ go to infinity, provided that $M \geq \frac{2N}{\varepsilon} \sup_{x \in \mathcal{S}_2}\mu( x)$.
\end{lemme}

\begin{pf}
Note that $M \geq \frac{2N}{\varepsilon} \sup_{x \in \mathcal{S}_2}\mu( x)$ implies
$(M+N) (\mu( y)-2 \epsilon) \le M (\mu(y)$ $-\,\epsilon)$ and use Lemma~\ref{lesboites} on the time needed to cross a box:
\[
p(M,N,\epsilon) \leq 2d \sup_{k \in \Zd} \sup_{y \in \mathcal{S}_2} \sup_{1
  \le m \le d} \P_p \bigl(t^m_{(y)}(k,\alpha r) \leq \bigl(\mu(
  y)-\epsilon\bigr)\alpha_m r\bigr).
  \]
\upqed\end{pf}

\section[Large deviations for the set of wet vertices: Proof of Theorem~1.7]{Large deviations for the set of wet vertices: Proof of Theorem~\protect\ref{dev-AS}}
\label{forme}
We can now prove Theorem~\ref{dev-AS}, which follows quite
naturally from the uniform estimates in Theorem~\ref{devsup}
and Theorem~\ref{devinf}.

Let $p>p_c(d)$ and $\epsilon>0$.
Let us note first that, for every $t>0$,
\begin{eqnarray*}
\Pcond_p  \biggl(
\mathcal{D} \biggl(\frac{B_t}t,\mathcal{B}_{\mu}(0,1) \biggr) \ge \epsilon
 \biggr) & \le &
\Pcond_p  \biggl(
\frac{B_t}t \not \subset \mathcal{B}_{\mu}(0,1+\epsilon)
 \biggr) \\
&&{} +\Pcond_p  \biggl(
\mathcal{B}_{\mu}(0,1) \not \subset \frac{B_t}t  +\mathcal{B}_{\mu}(0,\epsilon)
 \biggr).
\end{eqnarray*}
Let us now estimate each term separately.

\textit{Step }1. In the first term, we estimate the probability that the random set $B_t$ grows too fast, which corresponds to the existence of a point $x$ whose distance $D(0,x)$ from the origin is shorter than expected.
Thus
\begin{eqnarray*}
&&\Pcond_p  \biggl(
\frac{B_t}t \not \subset \mathcal{B}_{\mu}(0,1+\epsilon)
 \biggr)\\
&&\qquad =
\Pcond_p  \biggl(\exists x \in \Zd \, D(0,x)\le t, \,  \mu \biggl( \frac{x}{t} \biggr) >1 +\epsilon  \biggr) \\
&&\qquad =
\Pcond_p  \biggl(\exists x \in \Zd \, D(0,x)< \frac1{1+\epsilon} \mu(x), \,  \mu ( x ) >(1 +\epsilon)t  \biggr) \\
&&\qquad \le
\sum_{k={(1+\epsilon)t}/{\mu(e_1)}}^\infty \sum_{\|x\|_1=k}
\Pcond_p \biggl (
D(0,x)< \frac1{1+\epsilon} \mu(x)  \biggr)\\
&&\qquad \le
\sum_{k={(1+\epsilon)t}/{\mu(e_1)}}^\infty \sum_{\|x\|_1=k} Ae^{-B\|x\|_1}  \qquad     \mbox{with Theorem~\ref{devinf}}\\
&&\qquad \le  \sum_{k={(1+\epsilon)t}/{\mu(e_1)}}^\infty C k^{d-1} Ae^{-Bk} \le  A't^d e^{-Bt}.
\end{eqnarray*}

\textit{Step }2. In the second term, we estimate two types of discrepancies
between $B_t/t$ and $\mathcal B_\mu$: on the one hand the probability that
the random set $B_t$ grows too slowly, which corresponds to the existence
of a point $x$ whose distance $D(0,x)$ from the origin is larger than
expected and on the other hand, the probability that the random set $B_t$ contains
abnormally large holes. By definition,
\begin{eqnarray*}
&& \Pcond_p  \biggl(
\mathcal{B}_{\mu}(0,1) \not \subset \frac{B_t}t  +\mathcal{B}_{\mu}(0,\epsilon)
 \biggr)\\
&&\qquad =   \Pcond_p  \bigl(
\exists x \in \Rd \, \mu(x) \le t, \,  \mathcal{B}_{\mu}(x,\epsilon t)\cap B_t=\varnothing  \bigr).
\end{eqnarray*}
Note that, as soon as $t$ is large enough, one has
\begin{eqnarray*}
&&\forall x \in \mathcal{B}_{\mu}(0,t) \mbox{ } \exists y \in \Zd \mbox{ such that } \\
&&\qquad\frac{3\epsilon }4 t \le \mu(y) \le  \biggl( 1 -\frac{3\epsilon}4  \biggr) t
\quad\mbox{and}\quad
\mathcal{B}_{\mu} \biggl ( y, \frac{\epsilon t}8  \biggr) \subset
\mathcal{B}_{\mu}  ( x, \epsilon t  ),
\end{eqnarray*}
and so we obtain
\begin{eqnarray*}
&&\Pcond_p  \biggl(
\mathcal{B}_{\mu}(0,1) \not \subset \frac{B_t}t  +\mathcal{B}_{\mu}(0,\epsilon)
 \biggr)\\
&&\qquad \le   \Pcond_p \biggl (
\exists y \in \Zd \, \frac{3\epsilon }4 t \le \mu(y) \le  \biggl( 1
  -\frac{3\epsilon}4  \biggr) t ,    \mathcal{B}_{\mu} \biggl( y, \frac{\epsilon t}8  \biggr) \cap B_t=\varnothing  \biggr)\\
&  &\qquad\le  \sum_{({3\epsilon }/4) t \le \mu(y) \le  ( 1 -{3\epsilon}/4  ) t} \Pcond_p  \biggl(\mathcal{B}_{\mu} \biggl( y, \frac{\epsilon t}8  \biggr) \cap B_t=\varnothing  \biggr).
\end{eqnarray*}
If the event $\mathcal{B}_{\mu} ( y, \frac{\epsilon t}8  ) \cap B_t=\varnothing$ occurs, then either $\mathcal{B}_{\mu} ( y, \frac{\epsilon t}8  )$ contains no point of the infinite cluster, or it contains a point of the infinite cluster whose distance  from the origin is larger than $t$:
\begin{eqnarray*}
\Pcond_p  \biggl(\mathcal{B}_{\mu} \biggl( y, \frac{\epsilon t}8  \biggr) \cap B_t=\varnothing
\biggr)
& \le & \Pcond_p  \biggl(\mathcal{B}_{\mu}\biggl ( y, \frac{\epsilon t}8  \biggr) \cap C_\infty=\varnothing  \biggr) \\
&&{}+ \Pcond_p \biggl ( \exists z \in \mathcal{B}_{\mu} \biggl( y, \frac{\epsilon t}8  \biggr) \, t<D(0,z)<\infty  \biggr).
\end{eqnarray*}
By equation (\ref{amasinfini}), the first term is less than $A_3\exp(-\frac{B_3 \varepsilon t}{8 \mu(e_1)})$.
Note also that if $z \in \mathcal{B}_{\mu} ( y, \frac{\epsilon t}8  )$, then
$\mu(z)\le \mu(y)+\frac{\epsilon t}8 \le  ( 1-\frac{5\epsilon}8  )t$. Thus
\begin{eqnarray*}
&&\Pcond_p  \biggl( \exists z \in \mathcal{B}_{\mu} \biggl( y, \frac{\epsilon t}8  \biggr) \, t<D(0,z)<\infty  \biggr) \\
&&\qquad \le  \sum_{z \in \mathcal{B}_{\mu} ( y, {\epsilon t}/8  )} \Pcond_p  \biggl(
\biggl ( 1-\frac{5\epsilon}8  \biggr)^{-1} \mu(z)<D(0,z)<\infty  \biggr).
\end{eqnarray*}
For such $z$,  we have $\mu(z) \ge \mu(y)-\frac{\epsilon t}8 \ge \frac{5\epsilon}8 t$,
and by Theorem~\ref{devsup}, there exist two positive absolute constants $A,B$ such that
\[
\Pcond_p  \biggl( \exists z \in \mathcal{B}_{\mu} \biggl( y, \frac{\epsilon t}8
   \biggr) \,  t<D(0,z)<\infty  \biggr)
 \le  A(\epsilon t)^d e^{-B\epsilon t}.
\]
And we finally obtain
\begin{eqnarray*}
&&\Pcond_p  \biggl(
\mathcal{B}_{\mu}(0,1) \not \subset \frac{B_t}t  +\mathcal{B}_{\mu}(0,\epsilon)
 \biggr)\\
&&\qquad \le   \sum_{({3\epsilon }/4) t \le \mu(y) \le  ( 1 -{3\epsilon}/4  ) t} A_3\exp  \biggl(-\frac{B_3 \varepsilon t}{8 \mu(e_1)}  \biggr)+A(\epsilon t)^d e^{-B\epsilon t}\\
&&\qquad \le   Ct^d  \biggl( A_3\exp  \biggl(-\frac{B_3 \varepsilon t}{8 \mu(e_1)}  \biggr)+A(\epsilon t)^d e^{-B\epsilon t} \biggr),
\end{eqnarray*}
which ends the proof of the theorem.

\printaddresses

\end{document}